\documentclass[notitlepage, 11pt]{article}
\usepackage{graphicx}

\usepackage{bbm}
\usepackage{color}
\usepackage{latexsym}
\usepackage{amsmath, amssymb, amsfonts}
\usepackage[toc]{appendix}
\numberwithin{equation}{section}

\newtheorem{lemma}{Lemma}
\newtheorem{theorem}{Theorem}

\newtheorem{definition}{Definition}

\newtheorem{proposition}{Proposition}

\newtheorem{remark}{Remark}

\def\thelemma{\arabic{section}.\arabic{lemma}}

\newcommand{\beginsec}{
\setcounter{lemma}{0}
\setcounter{theorem}{0}
\setcounter{corollary}{0}
\setcounter{definition}{0}
\setcounter{example}{0}
\setcounter{proposition}{0}
\setcounter{condition}{0}
\setcounter{assumption}{0}
\setcounter{conjecture}{0}
\setcounter{problem}{0}
\setcounter{remark}{0}
}

\newcommand{\noi}{\noindent}

\newcommand{\E}{\mathbb{E}}
\newcommand{\R}{\mathbb{R}}

\newcommand{\N}{\mathbb{N}}

\newcommand{\la}{\lambda}
\newcommand{\sig}{\sigma}

\newcommand{\eps}{\varepsilon}

\newcommand{\Barr}{{\beta_\eps}}

\newcommand{\vr}{\varrho}
\newcommand{\al}{\alpha}
\newcommand{\gam}{\gamma}

\newcommand{\del}{\delta}
\newcommand{\om}{\omega}

\newcommand{\Gam}{\mathnormal{\Gamma}}

\newcommand{\X}{\mathnormal{\Xi}}

\newcommand{\Om}{\mathnormal{\Omega}}

\newcommand{\Q}{{\mathbb Q}}

\newcommand{\PP}{{\mathbb P}}

\newcommand{\calA}{{\cal A}}

\newcommand{\calC}{{\cal C}}
\newcommand{\calD}{{\cal D}}

\newcommand{\calF}{{\cal F}}
\newcommand{\calG}{{\cal G}}
\newcommand{\calH}{{\cal H}}
\newcommand{\calI}{{\cal I}}

\newcommand{\calL}{{\cal L}}

\newcommand{\calN}{{\cal N}}
\newcommand{\calO}{{\cal O}}

\newcommand{\calQ}{{\cal Q}}

\newcommand{\calX}{{\cal X}}

\newcommand{\skp}{\vspace{\baselineskip}}

\newcommand{\w}{\wedge}

\newcommand{\pl}{\partial}

\newcommand{\To}{\Rightarrow}

\newcommand{\dist}{{\rm dist}}

\newcommand\iy{\infty}
\newcommand{\osc}{\text{osc}}
\newcommand{\ds}{\displaystyle}

\newcommand{\qed}{\hfill $\Box$}

\newcommand{\limn}{\lim_{n\to\iy}}
\newcommand{\limm}{\lim_{m\to\iy}}
\newcommand{\limk}{\lim_{k\to\iy}}

\newcommand{\one}{\mathbbm{1}}
\newcommand{\aey}{a}

\newcommand{\uzero}{0}

\newcommand{\ue}{e}

\newcommand{\bB}{\boldsymbol{B}}

\newcommand{\bT}{\boldsymbol{T}}

\newcommand{\kaboom}{\hat\eps}

\title{Asymptotic analysis of a multiclass queueing control problem under heavy-traffic with model uncertainty\thanks{To appear in {\it Stochastic Systems}.}
}

\author{Asaf Cohen\thanks{Department of Statistics,
		University of Haifa,
		Haifa, 31905, Israel.
		Email:
		shloshim@gmail.com,
		web: https://sites.google.com/site/asafcohentau/. 
}}

\date{\today}

\begin{document}

\maketitle

\begin{abstract}
We study a multiclass M/M/1 queueing control problem with finite buffers under heavy-traffic where the decision maker is uncertain about the rates of arrivals and service
of the system and by scheduling and admission/rejection decisions acts to minimize a discounted cost that accounts for the uncertainty. The
main result is the asymptotic optimality of a $c\mu$-type of policy derived via underlying
stochastic differential games studied in \cite{Cohen2017}. 
Under this policy, with high probability, rejections are not performed when the workload lies below some cut-off that depends on the ambiguity level. When the workload exceeds this cut-off, rejections are carried out and only from the buffer with the cheapest rejection cost weighted with the mean service rate
. The allocation part of the policy is the same for all the ambiguity levels. 
This is the first work to address a heavy-traffic queueing control
problem with model uncertainty. 

\skp
\noi

\skp

\noi{\bf Keywords:} multiclass M/M/1, model uncertainty, heavy-traffic, diffusion approximation, Brownian control problem, stochastic differential game, state dependent priorities, ambiguity aversion.

\end{abstract}

\section{Introduction}\label{sec1}
\beginsec
The model under consideration consists of a multiclass M/M/1 queueing model under diffusion-scaled heavy-traffic where the decision maker (DM) is uncertain about the parameters and  acts to optimize an overall cost that accounts for this uncertainty. The model consists of a server that at any time instant her effort is allocated by the DM to customers from several number of classes. Customers of each class are kept in a finite buffer. Apart from the scheduling control, upon arrival of a customer, the DM has to decide whether to reject it or to assign it to the buffer that corresponds to its class type. The DM has ambiguity about the rates of arrivals and the mean service times. The cost accounts for the ambiguity, the holding of customers in the buffers, and rejections of new arrivals.

\skp\noi

The problem without ambiguity, also referred to as the {\it risk-neutral problem}, was analyzed by Atar and Shifrin in \cite{ata-shi}, under the framework of G/G/1. Plambeck et al.~studied in \cite{PKH} a similar non-robust problem with time constraints instead of the finite buffers constraints. In these problems as well as in many other classical models of queueing control problems (QCPs), see e.g., \cite{bell-will-1,BG2006,BG2012} and the references therein, there is a fixed time-homogeneous random model; that is, the DM is certain about the evolution of the system, which, moreover does not change in time. Such an assumption is not realistic, and a robust analysis is desirable. Due to the complexity of real-world systems, lack of sufficient calibration, and inaccurate assumptions, one cannot precisely model the arrival and departure processes. In recent years there has been increasing interest in robust analysis of queueing systems. We consider uncertainty in the diffusion scale of the QCP 
in a way that is often referred to as {\it model uncertainty} or {\it Knightian uncertainty}, see e.g., \cite{maenhout2004robust,Hansen2006,han-sar,bay-zha} and in the context of queueing systems see \cite{Shanti,blanchet2014robust,MR3544795}, see also \cite{MR1762853} in a discrete time setup. This is not to be confused with the terminology {\it robust queueing theory}, which is often referred to an optimization based performance analysis, see e.g., \cite{ban-ber-you2015,whi-you2018}. 
Uncertainty in fluid models of queueing was studied e.g., in \cite{whitt2006staffing,bassamboo2010capacity,dup2003}.

\skp\noi
We consider a DM who is skeptical about the validity of an underlying model and takes into account a family of models in the following way. We assume that the DM has a {\it reference model} in mind, which, up to some degree, describes the situation she is facing. To model the uncertainty about the reference model, the DM takes into account other models and penalizes them based on their deviation from the reference model. More specifically, we consider the following cost function, which the DM aims to minimize. 
\begin{align}\notag
\sup_{{\Q}}  &\left\{\E^{ \Q}\Big[\int_0^\iy e^{-\varrho t}\left(\hat h\cdot  X(t)dt +\hat r\cdot  d R(t)\right) \Big] - \alpha L(\Q\|P)
 \right\},
\end{align}
where $I$ is the number of classes. The vectors $\hat h$ and $\hat r$ stand for the holding and the rejection costs, respectively. The $I$-dimensional processes $X$ and $R$ represent the queue lengths and the rejection processes, respectively. The supremum is taken over measures equivalent to $\PP$,
and the function $L$ is a discounted variant of the Kullback--Leibler divergence. It measures how far away $\Q$ from the reference measure $\PP$ and penalizes accordingly. We allow for different levels of model uncertainty for each of the arrival processes as well as for each of the service processes. The parameter $\alpha>0$ measures the level of ambiguity the DM has regarding the reference model. It can also be viewed as the Lagrangian of a performance optimization problem with an uncertainty set constraint as presented in \cite[(14)]{Shanti} and \cite[(2)]{MR3544795}.
Moreover, exploiting the Markovian structure of the problem, we argue that one may replace the supremum over measures by a supremum over (possibly time-dependent) rates, which are penalized according to their deviation from the reference rates.


\skp\noi
To tackle QCPs under heavy-traffic one often solves a limiting control problem associated with a Brownian motion, called {\it Brownian control problem} (BCP) and uses its solution to construct an asymptotically optimal policy in the QCP. This concept was first introduced by \cite{har1988}; for further reading on BCPs see e.g., \cite{bell-will-1,BG2006,BG2012} and the references therein. In our case, the cost function given above suggests that the BCP is in fact a stochastic game. The players in this game are the DM and the nature, which according to their goals are referred to as the {\it minimizer} and the {\it maximizer}, respectively. Interpreting the roles of the processes from the QCP to a {\it multidimensional stochastic differential game}, the minimizer controls the server's effort allocation and the admission/rejection, while the maximizer is free to choose a probability measure and is penalized in accordance to the deviation of the chosen measure from the reference measure. This game was analyzed in \cite{Cohen2017}, where it was also shown that a state-space collapse property holds. This is done by considering a {\it reduced stochastic differential game} (RSDG), which emerges from the workload process, and showing that both games share the same value and that given equilibrium in either one of the games one can construct an equilibrium in the other game. Further properties of the games, such as dependency on the ambiguity parameters are also given there.

\skp\noi
This paper is devoted to the connection between the QCP and the BCP (namely, the RSDG); we show that the value function of the QCP is approximately the value function of the RSDG and that the minimizer's optimal strategy from the multidimensional stochastic differential game leads to an asymptotically Markovian optimal policy in the QCP.  Roughly speaking, this strategy suggests that the DM should use a $c\mu$ type of rule and fill in the buffers in accordance to their holding costs without using rejections, unless a cut-off level of workload that depends on the ambiguity level has been reached and then to use rejections only from the buffer with the cheapest rejection cost, weighted with the mean service rate, until the workload level goes below the cut-off. More specifically, let $\mu_i$ be the mean service rate of class $i$ customers under the reference model and recall that $\hat h_i$ is the holding cost per customer of class $i$. The DM should prioritize the classes in the order of $\{\hat h_i\mu_i\}_{i=1}^I$, where 
the lowest priority is given to the class with the lowest $\hat h_i\mu_i$ among the classes for which the buffers are not `almost' full. As for the admission control part of the policy, whenever the workload level remains below the mentioned cut-off use rejections only if there is a new arrival to a full buffer; as is shown, the probability of such an event vanishes with the scaling parameter. If the workload level exceeds the cut-off, the DM rejects all incoming customers to the class with the lowest $\hat r_i\mu_i$. 
This policy (with different cut-off level) was shown to be asymptotically optimal in the risk-neutral setup and in a QCP with the same mechanism but with the moderate-deviation heavy-traffic regime (instead of the diffusion scaling) and a risk-sensitive cost criteria. The latter QCP models a situation of a `very' risk averse DM. The only difference between the policies in the three models is in the the position of the cut-off point. Specifically, the allocation policy is the same. This shows us the usefulness of the allocation policy as it is robust to ambiguity. The optimality of the same policy (with differences only in the cut-off level) in these three models is not obvious due to the existence of the maximizer, which leads to a non-stationary problem. For further reading about the moderate-deviation heavy-traffic regime see \cite{ata-bis,bis2014,ata-coh,ata-sub_cmu,ata-coh2017, ata-men} and in the context of our paper also the discussion in  \cite[p.~3]{Cohen2017}. The current paper does not aim to establish a rigorous relationship between the QCP under the Knightian uncertainty and the QCP under the moderate-deviation heavy-traffic regime. Saying this, we only mention that a common feature to this model and the moderate-deviation one is that in the limiting problem the change of measure can be translated into a linear (in the maximizer's control variable) change in the drift and a quadratic penalty.

\skp\noi
We now make some comments on proof techniques. The proof is divided into two parts: showing that the value function of the RSDG forms an asymptotic lower bound for the QCP and that the expected cost associated with the candidate policy is asymptotically bounded above by the value function of the RSDG. 
Also, recall that we consider a sequence of queueing systems. For the lower bound we assign the maximizer a policy that is driven from the equilibrium strategy in the RSDG and using it we show that for any sequence of strategies used by the minimizer, the expected cost is bounded below by the value of the RSDG. By its structure, the maximizer's strategy preserves the critical load of the system. The main technical difficulty in this part is that the sequence of strategies of the minimizer is arbitrary and due to the nature of the control in the QCP, which is replaced by a singular control in the BCP, compactness arguments do not apply here. In the risk-neutral case, where there is no maximizing player, Atar and Shifrin \cite{ata-shi} managed to bypass this issue by $\calC$-tightness arguments applied to the integrands of the relevant processes and later on taking the derivatives of the implied limiting processes. In our setup, the maximizer's strategy depends on the scaled queue lengths process and not on its integrand. Therefore, the $\calC$-tightness of the sequence of the scaled queue lengths processes is required. As a result, the integrand-derivative method cannot be applied in our case. We use the {\it time-stretching} method, which was introduced by Meyer and Zheng in \cite{MZ1984} and studied in the same framework by Kurtz in \cite{Kurtz1991}. In the context of stochastic control the method was first used by Kushner and Martins in \cite{Martins1990, Kushner1991} and was adopted in \cite{BG2006, bud-ros2006, bud-ros2007}. In Section \ref{sec421} we set up a random time transformation for any system such that the controls are Lipschitz continuous with Lipschitz constant $1$. Then in Lemma \ref{lem41} we apply $\calC$-tightness arguments for the sequence of the relevant time-stretched processes (including the scaled queue lengths process) and obtain limiting processes. Using an inverse time transformation in Section \ref{sec422} we go back to the original scale. Finally, in Section \ref{sec423} we connect between the costs associated with the time rescaled limiting processes and the value of the RSDG. 

\skp\noi 
To show that the expected cost associated with the candidate policy is asymptotically bounded above by the value function of the RSDG we do the following. We 
consider an arbitrary sequence of strategies for the maximizer.  This sequence of probability measures is not forced to satisfy the critical load condition. Therefore, at first step we show that it is `too costly' for the maximizer to have a big deviation from the reference measure and thus restrict the maximizer to strategies that in average do not deviate much from the reference measure (Proposition \ref{lem44}). Then, in Proposition \ref{lem47} we adapt a state-space collapse result from the risk-neutral case to ours and show that under the sequence of strategies chosen by the maximizer,
the underlying stochastic dynamics of the scaled queue lengths lie close to a specific path, called the {\it minimizing curve}. 
At this point one can use $\calC$-tightness arguments to conclude the convergence of the relevant dynamics including the holding and rejection cost parts. However, in order to estimate the change of measure penalty given by the Kullback--Leibler divergence one shall have another reduction and `truncate' the maximizer's strategies such that the critically load condition is almost surely preserved. Then we show that the relevant processes and all the cost components associated with the two levels of restrictions of the maximizer's strategies are close to each other (Proposition \ref{prop41}). This approximation relies also on the state-space collapse. Thus, the rest of the analysis is perform in the more convenient case, where the critically load condition holds and the expected cost is shown to be asymptotically bounded from above by the value function of the RSDG. For this, we reduce to one-dimensional dynamics and estimate from below the change of measure penalty and together with the convergence of the holding and rejections cost components, we conclude the upper bound. Section \ref{sec5} provides a future outlook and Section \ref{sec6} is dedicated to the proofs of the more technical and less innovative lemmas and propositions stated in Section \ref{sec4}.

\skp\noi
The paper is organized as follows. Section \ref{sec2} presents the model. 
Section \ref{sec3} collects a few results from \cite{Cohen2017} required for the proof. In Section \ref{sec4} we provide and prove the main result (Theorem \ref{thm41}), which states that the QCP’s value converges to that of the BCP from \cite{Cohen2017} and an asymptotically optimal policy is provided. 

\subsection{Notation}\label{sec11}
We use the following notation.
For $a,b\in\R$, $a\wedge b=\min\{a,b\}$ and $a\vee b=\max\{a,b\}$. For a positive integer $k$ and $c,d\in\R^k$, $c\cdot d$ denotes the usual
scalar product and $\|c\|=(c\cdot c)^{1/2}$. We denote $[0, \iy)$ by $\R_+$. The infimum of the empty set is taken to be $\iy$. For subintervals $I_1,I_2\subseteq\R$ and $m\in\{1,2\}$ we denote by $\calC(I_1,I_2)$, $\calC^m(I_1,I_2)$, and $\calD(I_1, I_2)$ the space of continuous functions [resp., functions with continuous derivatives of order $m$, functions that are right-continuous with finite left limits (RCLL)] mapping $I_1\to I_2$. 
The space $\calD(I_1,I_2)$ is endowed with the usual Skorokhod topology. For $T,\delta>0$ and
a function $f:\R_+\to\R^k$, $\|f\|_T=\sup_{t\in[0,T]}\|f(t)\|$,
while $\osc_T(f,\del)=\sup\{\|f(u)-f(t)\|:0\le u\le t\le (u+\del)\w T\}$. 
For any RCLL processes $X,Y$, the quadratic variation of $X$ is denoted by $[X]$ and the quadratic covariation of $X$ and $Y$ is denoted by $[X,Y]$.

\section{The queueing model}\label{sec2}
\beginsec
We consider a QCP with customers of $I$ different classes that arrive at the system to be serviced by a single server. Upon arrival the customers are queued in $I$ buffers with finite capacity in accordance to their class.  Processor sharing is allowed and the server may serve up to $I$ customers at a time, where two customers from the same class cannot be served simultaneously. We study the system under heavy-traffic. Hence, we consider a sequence of systems, indexed by the scaling  parameter $n\in\N$.

\subsection{The reference model}\label{sec21}
Usually in (non-statistical) stochastic control problems it is assumed that the probability measure is fully known. A model is fixed and the DM finds the optimal policy given the probability measure. In statistics-based approach the parameters are not fully known and the DM learns them. We consider an epistemic {\it reference probability space} that is driven by previous knowledge as in non-statistical stochastic control problems. But instead of working only with this probability measure, as detailed in Section \ref{sec22}, the DM considers a robust cost that takes into account a family of measures without statistically learning the true model.

For every $i\in[I]:=\{1,\ldots,I\}$ and $n\in\N$ we consider the probability spaces
$(\Omega^n_{1,i},\calG^n_{1,i}
\PP^n_{1,i})$ and $(\Omega^n_{2,i},\calG^n_{2,i},
\PP^n_{2,i})$ that respectively 
support a Poisson process $A^n_i$ with a given rate $\lambda^n_i$
and a Poisson process $S^n_i$ with rate $\mu^n_i$. The process $A^n_i$ counts the number of arrivals to the $i$-th buffer 
%
and $S^n_i(t)$ stands for the number of service completions of customers of class $i$ had service was given to class $i$ for $t$ units of time.
Denote $A^n=(A^n_i)_{i=1}^I$ and $S^n=(S^n_i)_{i=1}^I$.

As will be rigorously clarified in Section \ref{sec22}, we consider a different level of uncertainty about each of the $2I$ components of $(A^n,S^n)$, and thus set the following {\it reference probability space} that supports these processes
\begin{align}\notag
(\Omega^n,\calG^n,
\PP^n):=
\Big(\prod_{i=1}^I(\Omega^n_{1,i}\times\Omega^n_{2,i}),\otimes_{i=1}^I(\calG^n_{1,i}\otimes\calG^n_{2,i})
,\prod_{i=1}^I(\PP^n_{1,i}\times\PP^n_{2,i})\Big),
\end{align}
where $\otimes_{i=1}^I(\calG^n_{1,i}\otimes\calG^n_{2,i})=(\calG^n_{1,1}\otimes\calG^n_{2,1})\otimes\ldots\otimes(\calG^n_{1,I}\otimes\calG^n_{2,I})$.

By the structure of the probability space, under the measure $\PP^n$, the $2I$ processes $A^n_1, S^n_1, \ldots,$ $A^n_I, S^n_I$ are mutually independent. Moreover, the distribution of $A^n_i$ (resp., $S^n_i$) under $\PP^n$ is identical to its distribution under $\PP^n_{1,i}$ (resp., $\PP^n_{2,i}$). 

Let $U^n=(U^n_i)_{i=1}^I$ be an RCLL process taking values in $\mathbb{U}=\{x=(x_1,\ldots,x_I)\in[0,1]^I:\sum x_i\le1\}$, where its $i$-th component $U^n_i(t)$ represents
the fraction of effort the server dedicates to class $i$ (recall that process sharing is allowed). Then
\begin{align}\label{206}
T^n_i(t):=\int_0^t U^n_i(s)ds,\qquad t\in\R_+,
\end{align}
gives the cumulative effort in the interval $[0,t]$ the server dedicates to class $i$ customers and the number of class $i$ job completions by time $t$ is thus $S^n_i(T^n_i(t))$. This is a Cox process with rate $\mu^n_iU^n_i$.

We allow rejections of customers (only) upon arrival and a rejected customer will never return to the system. The number of rejections from class $i$ until time $t$ is denoted by $R^n_i(t)$.

Denote by $X^n_i(t)$ the number of class $i$ customers in the system at time $t$. Then, the balance equation is given by
\begin{equation}\label{204a}
X^n_i(t)=X^n_i(0)+A^n_i(t)-S^n_i(T^n_i(t))-R^n_i(t),\qquad t\in\R_+.
\end{equation}
For simplicity, we assume that $X^n_i(0)$ is deterministic.
We use the notation $L^n=(L^n_i)_{i=1}^I$ for $L^n=X^n,R^n,T^n$. Since $U^n$ is an RCLL process, and by construction, so are $A^n$ and $S^n$, we conclude that $X^n$ and $R^n$ are RCLL as well.
The capacity of buffer $i$ is denoted by $ \hat b^n_i$ and is nondecreasing with $n$.

We now introduce the diffusion scaling and the heavy-traffic condition. First, we assume that
\begin{align}\label{205}
\lambda^n_i:=\lambda_i n+\hat\lambda_in^{1/2}+o(n^{1/2}),\qquad
\mu^n_i:=\mu_in+\hat\mu_in^{1/2}+o(n^{1/2}),
\end{align}
for some fixed constants $\lambda_i,\mu_i\in(0,\iy)$ and $\hat\lambda_i,\hat\mu_i\in\R$. Moreover, the system is assumed to be \emph{critically loaded}, that is,
$\sum_{i=1}^I\rho_i=1$, where $\rho_i:=\lambda_i/\mu_i$, $i\in[I]$.
Also, set up $\hat b^n_i:=\hat b_in^{1/2}$ for some constant $\hat b_i\in(0,\iy)$, $i\in[I]$.

The process $(U^n,R^n)$ is regarded as a control in the $n$-th system and is now given rigorously.

\begin{definition}[admissible control for the decision maker, QCP]\label{def21}
An {\it admissible control for the minimizer} for any initial state $ X^n(0)$ is a process $( U^n,R^n)$ taking values in $\mathbb{U}\times \R^I_+$ that satisfies the following,

\skp\noi
(i) $( U^n,R^n)$ is adapted to the filtration ${\calG}^n_t={\calG}^n(t):=\sigma\{A^n_i(s),S^n_i(T^n_i(s)), i\in[I], s\le t\}$ and has RCLL sample paths;

\skp \noi
(ii) the processes $R^n$ is nondecreasing with jumps of sizes $1/\sqrt{n}$;

\skp\noi
(iii) for each $i\in[I]$ and $t\in\R_+$,
\begin{align}\notag
X^n_i(t)=0\quad\text{implies}\quad U^n_i(t)=0;
\end{align}

\skp\noi
(iv) the buffer constraints $X^n_i(t)\in[0,\hat b_i^n]$, $t\in\R_+$, $i\in[I]$, hold.
\end{definition}
The first condition expresses the fact that the DM makes her decision based on past observations. The second condition follows since rejections are accumulated. The third condition asserts that service cannot be given to an empty buffer. We denote the set of admissible controls for the DM in the $n$-th system by $\calA^n( X^n(0))$.

\subsection{The optimization problem with model uncertainty}\label{sec22}

Recall that  we consider a DM that is uncertain about the underlying reference probability measure $\PP^n$, or loosely speaking, she suspects that the rates/intensities $\{\la^n_i\}_{i=1}^I$ and $\{\mu^n_i\}_{i=1}^I$ may be unspecified. Therefore, instead of optimizing under a single probability space, she considers a reference probability measure $\PP^n$ and a set of candidate measures (provided in the sequel) and penalizes their deviation from $\PP^n$. The penalization is done by using a discounted variant of the Kullback--Leibler divergence. More explicitly, the QCP is set up as a stochastic game that models a type of worst case scenario. The players are: the DM
that chooses a policy that minimizes a cost 
and an adverse player also referred to as the {\it nature} or {\it maximizer},
who has access to the policy chosen by the minimizer and to the history.
The nature is penalized for deviating from the reference model.

We are interested in a cost that accounts for the scaled queue lengths and rejections, in addition to the uncertainty about the model. Denote the scaled headcount process and the scaled rejection count by
\begin{align}\notag
&  \hat X^n(t):=n^{-1/2}X^n(t)\qquad\text{and} \qquad \hat R^n(t):=n^{-1/2}R^n(t),\qquad t\in\R_+.
\end{align}
Fix a discount factor $\varrho>0$, vectors of holding and rejection costs $\hat h,\hat r\in(0,\iy)^I$, and ambiguity parameters $\kappa:=(\kappa_{1,i},\kappa_{2,i})_{i=1}^I\in(0,\iy)^{2I}$. The DM chooses an admissible control $(U^n,R^n)$ in an attempt to minimize the following robust cost\footnote{Notice that this is a stochastic game, where the minimizer chooses a strategy and the maximizer, who knows the minimizer's choice, responses by choosing a worst case scenario. Also, we do not consider a statistical inference learning problem} 
\begin{align}\notag
\sup_{\hat{\Q}^n\in\hat{ \cal{Q}}^n( X^n(0))} J^n( X^n(0), U^n, R^n,\hat\Q^n;\kappa),
\end{align}
where the set of candidate measures $\hat\calQ^n( X^n(0))$ consists of all the product measures $\hat\Q^n=\prod_{i=1}^I(\hat\Q^n_{1,i}\times\hat\Q^n_{2,i})$ that satisfy some integrability condition (the exact definition is provided below), and 
\begin{align}\label{212}
 & J^n( X^n(0),U^n,R^n,\hat\Q^n;\kappa):=\\\notag
 &\quad\E^{\hat \Q^n}\Big[\int_0^\iy e^{-\varrho t}\left(\hat h\cdot\hat  X^n(t)dt +\hat r\cdot  d\hat R^n(t)\right) \Big]
 -\sum_{i=1}^I\frac{1}{\kappa_{1,i}}L^\varrho_1(\hat \Q^n_{1,i}\|\PP^n_{1,i})
 -\sum_{i=1}^I\frac{1}{\kappa_{2,i}}L^\varrho_2(\hat \Q^n_{2,i}\|\PP^n_{2,i}),
\end{align}
and
\begin{align}\label{211}
\begin{split}
L^\varrho_1(\hat \Q^n_{1,i}\|\PP^n_{1,i})&:=
\E^{\hat \Q^n_{1,i}}\left[\int_0^\iy\varrho e^{-\varrho t}\log\frac{d\hat \Q^n_{1,i}(t)}{d\PP^n_{1,i}(t)}dt\right],\\
L^\varrho_2(\hat \Q^n_{2,i}\|\PP^n_{2,i})&:=
\E^{\hat \Q^n_{2,i}}\left[\int_0^\iy\varrho e^{-\varrho t}\log\frac{d\hat \Q^n_{2,i}(t)}{d\PP^n_{2,i}(t)}dT^n_i(t)\right].
\end{split}
\end{align}
The first component in the cost function includes the holding and rejection costs under the measure $\hat \Q^n$. The last two sums in \eqref{212} are referred to as the {\it change of measure penalties}. While the second penalty depends on the control of the DM through the integrator $dT^n_i$, the same asymptotic results hold if we assume a penalty with integration $dt$. The main reasons for choosing this form is because it simplifies the analysis.

When $\kappa_{j,i}$ is `small' (resp., `big') we say that there is a weak (resp., strong) ambiguity about the rates of the processes $A^n_i$ and $S^n_i(T^n_i):=S^n_i(T^n_i(\cdot))$. 
The idea is that for small $\kappa_{j,i}$'s there is a big punishment per unit of deviation from the reference measure and therefore, the maximizer would prefer to choose measures
$\hat\Q^n_{j,i}$ close to $\PP^n_{j,i}$ (in the divergence sense), and as a consequence also the relevant expectations. 
However, one needs to make sure that the total punishment given by $\frac{1}{\kappa_{j,i}}L^\varrho_j(\hat \Q^n_{j,i}\|\PP^n_{j,i})$ is also small. In \cite[Theorems 5.1 and 5.2]{Cohen2017} we show that as the ambiguity parameters converge to zero, the stochastic differential games, which are provided in the same paper
, converge to the risk-neutral BCP studied in \cite{ata-shi}. Therefore, our problem indeed models ambiguity with respect to (w.r.t.) the risk-neutral model. 

The set of candidate measures $\hat\calQ^n( X^n(0))$ consists of all the product measures $\hat\Q^n=\prod_{i=1}^I(\hat\Q^n_{1,i}\times\hat\Q^n_{2,i})$ that
for every $i\in[I]$ and $ t\in\R_+$ the Radon--Nikodym derivatives satisfy
\begin{align}\label{216}
\begin{split}
\frac{d\hat \Q^n_{1,i}(t)}{d\PP^n_{1,i}(t)}&=\exp\Big\{\int_0^t\log\left(\frac{\psi^n_{1,i}(s)}{\la^n_i}\right) d A^n_i(s)-\int_0^t (\psi^n_{1,i}(s)-\la^n_i)ds\Big\}, \\
\frac{d\hat \Q^n_{2,i}(t)}{d\PP^n_{2,i}(t)}&=\exp\Big\{\int_0^t\log\left(\frac{\psi^n_{2,i}(s)}{\mu^n_i}\right) d S^n_i(T^n_i(s))-\int_0^t (\psi^n_{2,i}(s)-\mu^n_i)dT^n_i(s)\Big\},
\end{split}
\end{align}
for $\{\calG^n_t\}$-predictable measurable and positive processes $\psi^n_{j,i}$, $j\in\{1,2\}$, satisfying $\int_0^t\psi^n_{j,i}(s)ds<\iy$ $\PP$-almost surely (a.s.), for every $t\in\R_+$. 
These conditions assure, first, that the right-hand sides in \eqref{216} are $\PP^n_{j,i}$-martingales, $j\in\{1,2\}$ and second, that under the measure $\hat\Q^n_{1,i}$ (resp., $\hat\Q^n_{2,i}$), $A^n_i$ (resp., $S^n_i(T^n_i)$) is a counting process with intensity $\psi^n_{1,i}$ (resp., $\psi^n_{2,i}U^n_i$). Notice also that under the measures $\hat\Q^n_{j,i}$, $j\in\{1,2\}$, the critically load condition might be violated since we do not restrict the intensities $\{\hat\psi^n_{j,i}\}_{j,i,n}$ in such a way. However, as we detail in Section \ref{sec44} below 
such changes of measures are `too costly' and will be avoided by the maximizer so that, {\bf as a consequence}, `in average' the critically load condition is preserved. 

Finally, the value function is thus given by
\begin{align}\notag
 V^n( x;\kappa):= \;&
\inf_{(U^n,R^n)\in\calA^n( x)}\;
\sup_{\hat{\Q}^n\in\hat{ \cal{Q}}^n( x)} J^n( x, U^n, R^n,\hat\Q^n;\kappa).
\end{align}

\begin{remark}\label{rem_21}
(i) Notice that the right-hand side (r.h.s.) of \eqref{216} are martingales, and therefore, there exist probability measures $\hat\Q^n_{j,i}$, $j\in\{1,2\}$, such that $\hat\Q^n_{j,i}|_{\calG^n_t}$ satisfies \eqref{216} for all $t\in\R_+$, see \cite[Lemma 4.2]{Stroock1987}.

(ii) From the presentations in \eqref{216} it follows that one may think of the maximizer as choosing rates rather than choosing a measure. The reason is that the rates $(\la^n_i,\mu^n_i)$ under the reference measure $\PP$ are replaced by the rates $(\psi^n_{1,i}(\cdot),\psi^n_{2,i}(\cdot))$, which are allowed to be time-dependent. Including time-dependent change of measures/rates is more realistic, since in real world systems one has no guarantee that the rates remain fixed over time, and moreover, we show that a $c\mu$-type of policy remains optimal under this general structure.

%
%

\end{remark}
%

\section{The BCP}\label{sec3}
\beginsec
\subsection{The reduced stochastic differential game (RSDG)}\label{sec32}
In \cite{Cohen2017} we provided two equivalent games: an $I$-dimensional game that captures the dynamics of the buffers and a one-dimensional game, whose dynamics approximate the workload process. For the forthcoming analysis we will use only the latter. This section is devoted to the description and the relevant properties of this game. 

To derive the game we introduce some scaled processes and constants in addition to $\hat X^n$ and $\hat R^n$ introduced earlier. Set for every $i\in[I]$ and $t\in\R_+$
\begin{align}\label{207}
\begin{split}&
\hat{A}^n_i(t):=n^{-1/2}(A^n_i(t)-\la^n_i t),\qquad
 \hat{S}^n_i(t):=n^{-1/2}(S^n_i(t)-\mu^n_i t),\\
 &
\hat Y^n_i(t):=\mu^n_in^{-1/2}(\rho_i t-T^n_i(t)),\qquad 
\hat m^n_i:=n^{-1/2}(\lambda^n_i-\rho_i\mu^n_i).
\end{split}
\end{align}
We use the notation
$ \hat X^n=(\hat X^n_i)_{i=1}^I$ and similarly for $\hat R^n,\hat A^n,\hat S^n,\hat Y^n,$ and $\hat m^n$. Denote also $\hat S^n(T^n)=(\hat S^n_i(T^n_i(\cdot)))_{i=1}^I$.
The scaled version of \eqref{204a} is given by,
\begin{equation} \label{new18}
 \hat{X}^n(t)=\hat{X}^n(0) + \hat m^nt
 +\hat{A}^n(t)-\hat{S}^n(T^n(t))
  +\hat  Y^n(t)- \hat{R}^n(t),\qquad t\in\R_+.
\end{equation}
An admissible policy satisfies
\begin{align}\label{208}
\hat X^n(t)\in\calX:=\prod_{i=1}^I[0, \hat b_i],\qquad t\in\R_+,\;\; \PP^n\text{-a.s.}
\end{align}
Multiplying both sides of \eqref{new18} by the workload vector
\begin{align}\label{thetan}
\theta^n:=(n/\mu^n_1,\ldots,n/\mu^n_I),
\end{align}
we get
\begin{equation}\notag
\theta^n\cdot \hat{X}^n(t)=\theta^n\cdot\hat{X}^n(0) +\theta^n\cdot \hat m^nt
 +\theta^n\cdot\left(\hat{A}^n(t)-\hat{S}^n(T^n(t))\right)
  +\theta^n\cdot\hat  Y^n(t)- \theta^n\cdot\hat{R}^n(t).
\end{equation}
Under the reference measure $\PP^n$, $\{(\hat A^n,\hat S^n)\}_n$ weakly converges a $2I$-dimensional $(0,\tilde\sigma)$-Brownian motion, where $\tilde\sigma:=\text{Diag}\;(\la_1^{1/2} ,\ldots,\la_I^{1/2},\mu_1^{1/2} ,\ldots,\mu_I^{1/2}).$ 
As we show rigorously in the proof of Lemma \ref{lem45}, $\hat Y^n$ is of order one as $n\to\iy$. Hence its definition implies that $T^n(t)\to(\rho_1,\ldots,\rho_I)t$, $t\in\R_+$, and therefore, under $\PP^n$, $\{\hat A^n-\hat S^n(T^n)\}_n$ weakly converges to an $I$-dimensional $(0,\hat\sigma)$-Brownian motion, where 
\begin{align}\notag
\hat\sigma=(\hat\sigma_{ij}):=\text{Diag}\left((2\la_1)^{1/2} ,\ldots,(2\la_I)^{1/2}\right).
\end{align} 
Hence, $\theta^n\cdot\left(\hat{A}^n(t)-\hat{S}^n(T^n(t))\right)$ can be approximated by a one-dimensional Brownian motion. This last two approximations suggest that the $2I$ ambiguity parameters $\{\kappa_{i,j}\}_{i,j}$ collapse into $I$ parameters $\{\eps_i\}_i$, which in turn are folded into a single parameter $\eps$ as follows:
\begin{align}\label{300b}
\kaboom_i:=\frac{1}{2}(\kappa_{1,i}+\kappa_{2,i}),\qquad i\in[I],\qquad\text{and}\qquad \eps:=\frac{1}{\sigma^2}\sum_{i=1}^I(\theta\hat\sigma)_i^2\kaboom_i.
\end{align}

 To introduce the game, we need the following notation. Set the vectors
\begin{align}\label{theta}
\theta:=(\mu_1^{-1},\ldots,\mu_I^{-1}),\qquad \hat b=(\hat b_1,\ldots, \hat b_I),
\end{align}
and the scalars
\begin{align}\label{310b}
x_0&:=\theta\cdot \hat x_0,\quad m:=\theta\cdot\hat  m,\quad\sigma:=\|\theta\hat\sigma\|,\quad b:=
\theta\cdot \hat b.
\end{align}
We are now ready to define the {\it reduced stochastic differential game} (RSDG).

\begin{definition} [admissible controls, RSDG]\label{def32}
An {\it admissible control for the minimizer} for any initial state $x_0\in[0,b]$ is a filtered probability space  $(\Omega,\calF,\{\calF_t\},\PP)$ that supports a one-dimensional standard Brownian motion $ B$ and a process $ ( Y, R)$ taking values in $\R_+^2$ with RCLL sample paths, both adapted to the filtration $\{\calF_t\}$ and satisfy the following properties:

\noi
(i) for every $0\le s<t$, $ B(t)- B(s)$ is independent of $\calF_s$ under $\PP$;

\noi
(ii) $ Y$ and $ R$ are nonnegative and nondecreasing;

\noi
(iii) The controlled process satisfies $
 X(t)= x_0+ mt+\sigma  B(t)+ Y(t)- R(t),
$ and
$ X(t)\in [0,b],\; t\in\R_+,\;\PP\text{-a.s.}
$

An {\it admissible control for the maximizer} is a measure $ \Q$ defined on $(\Omega,\calF,\{\calF_t\})$ such that
\begin{align}\label{315}
\frac{d \Q(t)}{d\PP(t)}=\exp\Big\{\int_0^t \psi(s)dB(s)-\frac{1}{2}\int_0^t  \psi^2(s)ds\Big\},\quad t\in\R_+,
\end{align}
for an $\{\calF_t\}$-progressively measurable process $\psi$ satisfying
\begin{align}\label{316}
&\E^{\PP}\Big[\int_0^\iy e^{-\varrho s} \psi^2(s)ds\Big]<\iy\quad\text{and}\quad\E^{\PP}\Big[e^{\frac{1}{2}\int_0^t \psi^2(s)ds}\Big]<\iy\quad\text{ for every $t\in\R_+$.}
\end{align}
\end{definition}

Denote by $ \calA( x_0)$ (resp., $ \calQ( x_0)$) the set of all admissible controls for the minimizer (resp., maximizer), given the initial condition $ x_0$. In \cite[(2.29)]{Cohen2017} it is argued that the controlled process can alternatively be written as 
\begin{align}\label{asaf001}
	X(t)= x_0+ mt+\int_0^t\sigma\psi(s)ds+ \sigma  B^{\Q}(t)+ Y(t)- R(t),\quad t\in\R_+,
\end{align}
where $ B^{ \Q}(t):= B(t)-\int_0^t\psi(s)ds$, $t\in\R_+$, is an $\{\calF_t\}$-one-dimensional standard Brownian motion under $\Q$. This form serves us in the asymptotic analysis.

The cost associated with the initial condition $x_0$ and the controls $( Y, R)$ and $ \Q$ is given by
\begin{align}\notag
J( x_0, Y, R, \Q;\eps):=&
\E^{ \Q}\Big[\int_0^\iy e^{-\varrho t}\left( h(  X(t))dt + rd R(t)\right) \Big]-\frac{1}{\eps}L^\varrho( \Q\|\PP),
\end{align}
where
\begin{align}\label{319a}
&h(x):=\min\{\hat h\cdot\xi : \xi\in\calX,\; \theta\cdot\xi=x\},\\\label{319b}
&r:=\min\{\hat r\cdot q : q\in\R^I_+,\; \theta\cdot q=1\},
\end{align}
and 
\begin{align}\label{321}
L^\varrho( \Q\|\PP):=\E^{\hat \Q}\left[\int_0^\iy\varrho e^{-\varrho t}\log\frac{d\hat \Q(t)}{d\PP(t)}dt\right] 
=\E^{ \Q}\Big[\int_0^\iy e^{-\varrho t} \hat\psi^2(t)dt \Big],
\end{align}
where the second inequality is given in \cite[(2.23)]{Cohen2017}. The latter form will turn out to be useful in the approximation procedure, see \eqref{newnew21} and \eqref{newnew20} below. Notice that the scalar $\eps$ folds all the ambiguity parameters $\{\kappa_{i,j}\}_{i,j}$ and all the $2I$ penalties from \eqref{212} are collapsing into one penalty. 
By the convexity of $\calX$ it follows that $h$ is convex. In fact, $h$ is piecewise linear and Lipschitz continuous. Moreover, $h(x)\ge 0$ for $x\ge 0$ and equality holds if and only if $x=0$. Therefore, $h$ is strictly increasing. In \cite[page 568]{ata-shi} it is shown that there is $i^*\in[I]$ such that,
\begin{align}\label{320}
r=r_{i^*}\mu_{i^*}:=\min\{\hat r_i\mu_i : i\in[I]\}.
\end{align}
The index $i^*$ stands for the class with the smallest rejection cost, weighted with the mean service rate. In fact, under the candidate asymptotically nearly optimal policy, presented in Section \ref{sec41}, unless a buffer is overloaded, an event that happens with vanishing probability, rejections are performed only from class $i^*$.

The value function is given by
\begin{align}\label{322}
 V( x_0;\eps)=\inf_{( Y, R)\in \calA( x_0)}\;\sup_{ \Q\in \calQ( x_0)}\;J( x_0, Y, R, \Q;\eps).
\end{align}

\subsection{Properties of the game}\label{sec33}
The RSDG admits a simple optimal strategy for the minimizer that enforces the workload to stay in a specific interval of the form $[0,\beta]$ with minimal effort. To rigorously define such a strategy we make use of the {\it Skorokhod map on an interval}. Fix $\beta>0$. For any $\eta\in\calD(\R_+,\R)$ there exists a unique triplet of functions $(\chi,\zeta_1,\zeta_2)\in\calD(\R_+,\R^3)$ that satisfies the following properties:

\noi
(i) for every $t\in\R_+$, $\chi(t)=\eta(t)+\zeta_1(t)-\zeta_2(t)$;

\noi
(ii) $\zeta_1$ and $\zeta_2$ are nondecreasing, $\zeta_1(0-)=\zeta_2(0-)=0$, and
\begin{align}\notag
\int_0^\iy
\one_{(0,\beta]}(\chi(t))d\zeta_1(t)=\int_0^\iy
\one_{[0,\beta)}(\chi(t))d\zeta_2(t)=0,
\end{align}
where $\one_F(x)=1$ if $x$ belongs to the set $F$ and 0 otherwise. We denote by $\Gamma_{[0,\beta]}(\eta)=(\Gamma_{[0,\beta]}^1,\Gamma_{[0,\beta]}^2,\Gamma_{[0,\beta]}^3)(\eta)=(\chi,\zeta_1,\zeta_2)$.
See \cite{Kruk2007} for existence and uniqueness of the solution,
and continuity and further properties of the map.
In particular, we have the following.
\begin{lemma}\label{lem_Skorokhod}
There exists a constant $c_S>0$ such that
for every $t>0$, $\beta,\tilde\beta>0$ and $\om,\tilde\om\in\calD(\R_+,\R)$,
\begin{equation}\notag
\|\Gam_{[0,\beta]}(\om)(t) - \Gam_{[0,\tilde\beta]}(\tilde\om)(t)\|_t
\le c_S(\|\om-\tilde\om\|_T+|\beta-\tilde\beta|)
.
\end{equation}
\end{lemma}

\begin{definition}\label{def_Skorokhod}
Fix $x_0,\beta\in[0,b]$. 
The strategy $(Y,R)$ is called 
a $\beta$-reflecting strategy if for every $\eta\in \calC(\R_+,\R)$ one has $(X,Y,R)(\eta)=\Gamma_{[0,\beta]}(\eta)$, with $X(0)=x_0$. 
\end{definition}

The next proposition summarizes the necessary results for the analysis in the sequel.
\begin{proposition}[Theorems 3.1 and 4.1 in \cite{Cohen2017}]\label{prop31}
Fix $\eps\in(0,\iy)$. The value function $V(\cdot):=V(\cdot;\eps)$ is twice continuously differentiable and satisfies $0\le V'\le r$. Moreover, set the parameter
\begin{align}\label{440}
\beta_\eps=\inf\left\{x\in(0,b] : V'(x)=r\right\}\wedge b
\end{align}
and the measure 
 $\Q_V=\Q_{V(\cdot;\eps)}$ associated with $\psi_V(t):=
\eps\sigma V'(X(t))$ via \eqref{315}.
Then, the reflecting strategy $(Y_{\beta_\eps},R_{\beta_\eps})$ and the measure $\Q_V$ form an equilibrium in the sense that
\begin{align}\notag
&V(x_0)=
\sup_{\Q\in\calQ(x_0)}J(x_0,Y_{\beta_\eps},R_{\beta_\eps},\Q;\eps)=
J(x_0,Y_{\beta_\eps},R_{\beta_\eps},\Q_V;\eps)=
\inf_{(Y,R)\in\calA(x_0)}J(x_0,Y,R,\Q_V;\eps).
\end{align}
\end{proposition}
In \cite{Cohen2017} it is shown that the value function $V$ is continuous and decreasing in the ambiguity parameter $\eps$ and that the reflecting barrier $\beta_\eps$ is continuous in $\eps$. The meaning of the set of equations in the proposition above is that by using the reflecting strategy $(Y_{\beta_\eps},R_{\beta_\eps})$, the best the maximizer can do is choosing the measure $\Q_V$, and vice versa. The optimality of $\Q_V$ serves us in Section \ref{sec42}, where we equip the maximizer in the QCP with a measure driven by $\Q_V$. The reflecting strategy $(Y_{\beta_\eps},R_{\beta_\eps})$ is used to show the state-space collapse in Section \ref{sec432}.

\section{Asymptotic analysis}\label{sec4}
\beginsec
\subsection{Nearly optimal policy}\label{sec41}
In this section we present a candidate policy that is shown to be nearly optimal. 
The idea of the policy is to give the least priority to the buffer with the cheapest holding cost and let it fill up until it is {\it almost} full. If the workload increases, we keep the cheapest buffer full and start filling the next cheapest buffer by less prioritizing it, and so on. This is done with {\it almost} no rejections, unless the critical workload level $\Barr$ (see Proposition \ref{prop31}) is reached in which case rejections occur, and only from the buffer with the cheapest rejection cost. 

\subsubsection{The minimizing curve}
\begin{figure}
	\centering
	\includegraphics[width=0.5\textwidth]{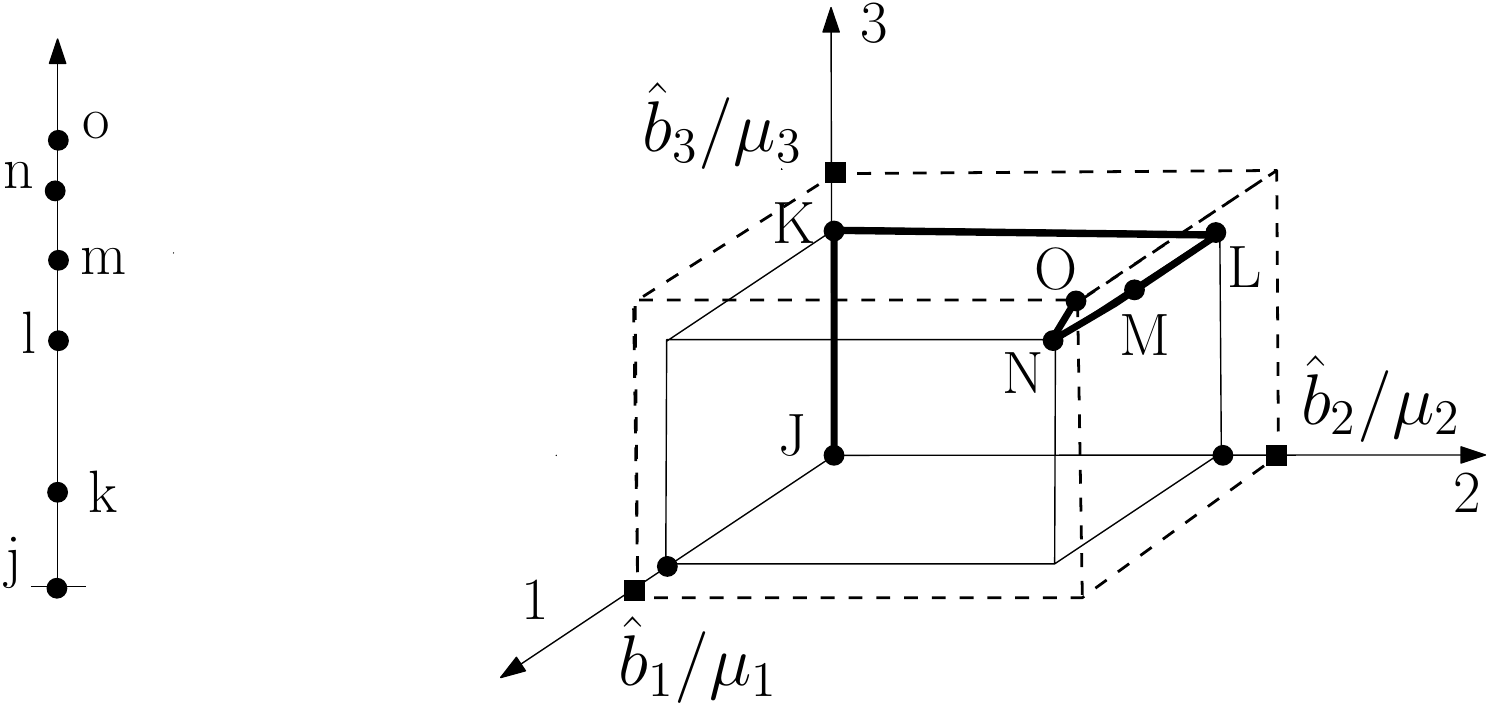}
	\caption{{\footnotesize The graphs refer to the case $I=3$, $\hat h=(1,5/2,3/2)$, $\mu=(3,1,3/2)$, and $(\hat a_1,\hat a_2,\hat a_3)=(4,7,6)$. The graph to the left stands for the workload levels. The curve of the function $\gamma^a$ is in bold in the graph to the right. The workload levels with the lower case letters are
			$j=0, k=\hat a_3/\mu_3=4, l=\hat a_3/\mu_3+\hat a_2/\mu_2=11, m= \hat a_3/\mu_3+\hat a_2/\mu_2+1=12$, and $n=\hat a_3/\mu_3+\hat a_2/\mu_2+\hat a_1/\mu_1=a=37/3$. They respectively correspond to the upper case letters: $J=(0,0,0), K=(0,0,\hat a_3/\mu_3)=(0,0,4), L=(0,\hat a_2/\mu_2,\hat a_3/\mu_3)=(0,7,4), M=(1,\hat a_2/\mu_2,\hat a_3/\mu_3)=(1,7,4)$, and $N=(\hat a_1/\mu_1,\hat a_2/\mu_2,\hat a_3/\mu_3)=(4/3,7,4)$.
			Pay attention that between $N$ and $O$, $\gamma^a$ is defined via a linear interpolation. If for example $m= \beta_\eps=a$, the scaled queues will fluctuate close to the bold curve between the points $J$ and $M$.
		}\label{fig_1}}
\end{figure}

Before presenting the candidate policy, we provide a {\it minimizing curve} that the dynamics in the QCP $\hat X^n$ asymptotically trace. Its structure stems from the equilibrium strategy in the multidimensional stochastic differential game provided in \cite{Cohen2017}. Due to the discrete nature of the QCP we take this curve to be bounded away from (yet close to) the boundary of $\Pi_{i=1}^I[0,\hat b_i/\mu_i]$. This curve is a key element in establishing the state-space collapse property. Without loss of generality assume that 
\begin{align}
\label{new4}
h_1\mu_1\ge h_2\mu_2\ge\dots\ge h_I\mu_I.
\end{align}
Fix $\delta_0>0$ and
let $\hat a=(\hat a_1,\ldots,\hat a_I)$ be given by $\hat a_i=\hat b_i-\delta_0$, $i\in[I]$, and 
$\aey =\Barr\w(\theta\cdot \hat a)<b$. Note that if $\delta_0$ is sufficiently small
then $\aey =\Barr$ (unless $\Barr=b$). The difference between the intercepts on axis $i=1,2,3$ in Figure \ref{fig_1} equals $\delta_0/\mu_i$. Hence, when $\delta_0\to 0$, the solid box converges to the dotted one.
We now rigorously define a function that maps workload values to states of the $I$-th dimensional system and then explain its structure. Set the function $\gamma^a:[0,b]\to\calX$ as follows. For $x\in[0,\theta\cdot\hat a)$, the variables $j=j(x)$ and $\upsilon=\upsilon(x)$ are determined
via
\begin{align}\label{eq2304}
x=\sum_{i=j+1}^I\theta_i\hat a_i+\theta_j\upsilon, \qquad j\in[I],\quad\upsilon\in[0,\hat a_j),
\end{align}
and let
\begin{align}\label{asaf003}
\gamma^a(x)=\sum_{i=j+1}^I\hat a_i\ue_i+\upsilon \ue_j.
\end{align}
On the interval $[\theta\cdot \hat a,b]$
define $\gam^a$ as the linear interpolation between the points $(\theta\cdot\hat a,\hat a)$
and $( b,\hat b)$. Notice that the buffers $j+1,\ldots, I$, which are the cheapest ones (holding cost-wise), are almost full (up to $\delta_0$), buffer $j$ is building up, and the buffers $1,\ldots, j$ are left empty.
An illustration of the curve $\gamma^a$ is provided in Figure \ref{fig_1}.

We now set up an holding cost that is associated with workload values through the curve $\gamma^a$. 
Recall the definition of $h$ from \eqref{319a} and let
\[
h^a(x):=\min\{\hat h\cdot\xi:\xi\in\calX, \theta\cdot\xi=x,
\xi_i\le\gamma^a_i(x),i\in[I]\}=\hat h\cdot\gamma^a(x),\qquad
x\in[0, \theta\cdot a].
\]
Also, set
\begin{align}\label{eq2307b}
\om_1(\delta_0)=\sup_{[0,\theta\cdot \hat a]}| h^a- h|.
\end{align}
By the choice of $a$ it is clear that $\om_1(0+)=0$.

\subsubsection{The candidate policy and the main theorem}
{\it Rejection policy:}
In case that a class-$i$ arrival occurs at a time $t$ when $\hat X^n_i(t-)+n^{-1/2}>\hat b_i$,
then it is rejected.
Such rejections are called \emph{forced rejections}.
Whenever $\theta\cdot \hat X^n\ge \aey $, all class-$i^*$ (see the paragraph preceding \eqref{320}) arrivals are rejected,
and these rejections are called \emph{overload rejections}. 
Apart from that, no rejections occur from any class.

\skp\noi{\it Service policy:}
For each $\hat x=(\hat x_1,\ldots,\hat x_I)\in\calX$ define the class of low priority
\[
\calL( \hat x)=\max\{i:\hat x_i<\hat a_i\},
\]
provided $\hat x_i<\hat a_i$ for some $i$, and set $\calL( \hat x)=I$ otherwise.
The complement set is the set of high priority classes:
\[
\calH(\hat  x)=[I]\setminus\calL( \hat x).
\]
When at least one class among $\calH( \hat x)$
is not empty, the class 
$\calL( \hat x)$ receives no service, and
all classes within $\calH( \hat x)$ that are not empty 
receive service at a fraction proportional to their
traffic intensities. Namely, denote $\calH^+(\hat  x)=\{i\in \calH( x):x_i>0\}$,
and define $\rho'( \hat x)\in\R^I$ as
\begin{align}\label{eq2310}
\rho'_i( \hat x)=\begin{cases}
0, & \text{if } \hat  x=\uzero,\\
\ds\frac{\rho_i\one_{\{i\in\calH^+( \hat x)\}}}{\sum_{k\in\calH^+( \hat x)}\rho_k}, &
\text{if } \calH^+( \hat x)\ne\emptyset,\\
e_I,& \text{if $\hat x_i=0$ for all $i<I$ and $\hat x_I>0$,}
\end{cases}
\end{align}
where recall that $e_I=(0,\ldots,0,1)\in\R^I$. Note that $\calH^+( \hat x)=\emptyset$ can only happen if $\hat x_i=0$ for all $i<I$,
which is covered by the first and last cases in the above display.
Then for each $t\in\R$,
\begin{align}\label{eq2311}
  U^n(t)=\rho'(\hat X^n(t)).
\end{align}
Note that when $\calH^+( \hat x)\ne\emptyset$,
\begin{align}\label{eq2312}
  \rho'_i( \hat x)>\rho_i\quad \text{for all } i\in \calH^+( \hat x).
\end{align}
That is, all prioritized classes receive
a fraction of effort strictly greater than the respective traffic intensity. 
Also note that $\sum_iU^n_i=1$ whenever $\hat X^n$ is nonzero.  
This is therefore a work conserving policy. 

To illustrate the relationship between the candidate policy and the minimizing curve we assume that the system starts empty. Then the low priority is given to Buffer $I$ (unless in times it is empty) and it fluctuates, while the other buffers get priority, hence are almost empty. This is the situation until Buffer $I$ exceeds the level $\hat a_I$. Then, Buffer $I$ receives priority and decreases rapidly until it goes below level $\hat a_I$ and now again is low prioritized. This way, while Buffer $I-1$ is not empty Buffers $I$ and $I-1$ exchange priorities so that Buffer $I$ remains around the level $\hat a_I$ and Buffer $I-1$ fluctuates. 

\begin{remark}
	Pay attention that the service policy does not depend on the ambiguity but rather only on the prioritization of the buffers according to the holding costs. The ambiguity, which is folded into the parameter $\eps$ plays a role in the rejection level $a$, which equals $\beta_\eps$, unless $\beta_\eps=b$, in which case $a$ is very close to $\beta_\eps$. As mentioned in the introduction, the candidate policy (again, with different rejection level) is also optimal under the moderate-deviation heavy traffic regime. This fact illustrates the robustness of the service policy, especially since our formulation allows for time-varying rates.	
\end{remark}

Recall the definition of $\theta^n$ given in \eqref{thetan}.
\begin{theorem}\label{thm41} Assuming that $x_0:=\limn\theta^n\cdot \hat X^n(0)$ exists, then,
 \begin{align}\label{new6}
  \lim_{n\to\iy}V^n( X^n(0);\kappa)= V(x_0;\eps).
  \end{align}
Moreover,
for every $n\in\N$, denote the policy constructed above by $(U^n(a),R^n(a))$. Then, there is a function $w_2:\R_+\to\R$, satisfying $\omega_2(0+)=0$ such that
\begin{align}\label{new5}
\limsup_{n\to\iy}\sup_{\hat{\Q}^n\in\hat{ \cal{Q}}^n( X^n(0))} J^n(X^n(0), U^n(a), R^n(a),\hat\Q^n;\kappa)\le V(x_0;\eps)+\omega_2(\Barr-a).
\end{align}
\end{theorem}
Recall that the parameter $a$ can be chosen to be arbitrary close to $\beta_\eps$, see the paragraph below \eqref{new4}. Therefore, by a diagonalisation argument, one can deduce an asymptotically optimal policy generated from $U^n(a)$.
The proof of the theorem takes place in the next two sections. In Section \ref{sec42} we show that the game's value function $V$ bounds from below the liminf of the QCP's value functions. That is,\begin{align}\label{new7}
\liminf_{n\to\iy}V^n(X^n(0);\kappa)\ge V(x_0;\eps).
\end{align}
In Section \ref{sec43} we prove \eqref{new5}. Together, we obtain \eqref{new6}.








\subsection{Proof of \eqref{new7}}\label{sec42}

%
The proof follows along these lines. We equip the maximizer in the QCP with a measure in such a way that its one-dimensional folding onto the workload scale is the equilibrium measure $\Q_V$ from Proposition \ref{prop31}. The DM is equipped with an arbitrary sequence of policies $\{(\hat Y^n,\hat R^n)\}_n$. Using weak convergence arguments and projecting the $I$-dimensional processes on the workload vector $\theta^n$ (see \ref{thetan}) we show that along every converging subsequence the projected dynamics converges to the dynamics given in \eqref{asaf001}, with $\psi=\psi_V$ (see Proposition \ref{prop31}), and establish the lower bound for the value functions.

The main difficulty of this argument is that since $\{(\hat Y^n,\hat R^n)\}_n$ is an arbitrary sequence of singular controls 
one cannot expect this sequence to be tight. Therefore, we use time stretching in order to prove the lower bound.\footnote{At this point, it is worth mentioning that Atar and Shifrin \cite{ata-shi}  managed to bypass this issue by arguing $\cal C$-tightness of the integrands of the relevant processes. Repeating the same arguments, one may show that the integrands of $\hat X^n_i$, $\hat Y^n_i$, and $\hat R^n_i$ are $\cal C$-tight. Since $V'$ is bounded, see Proposition \ref{prop31}, we get that the sequence $\left\{\int_0^t\left(\hat\psi^n_{1,i}(s)-\hat\psi^n_{2,i}(s)\right)ds\right\}_n$ is $\cal C$-tight. However, since 
we wish to obtain (after projecting on the workload vector) the specific integral given in \eqref{asaf002} below, we still need to argue tightness of $\hat X^n$.}
The idea is as follows. Take for example the process $\hat R^n$. In Section \ref{sec421} we stretch the time (using the same transformation for all the processes) and generate a process $\tilde R^n$ in such a way that $\{\tilde R^n\}$ is Lipschitz-continuous with Lipschitz constant 1 (over intervals of order 1) and therefore tight and converges to a process $\tilde R$. Then in Section \ref{sec422} we go back to the original scale by an inverse time transformation and get the process $\hat R$, which is used in Section \ref{sec423} to get the value function of the RSDG.

We start with setting up the maximizer's measure in the QCP.  
Recall the definition of $\theta$ from \eqref{theta}. For every $t\in\R_+$ and $i\in[I]$, set
\begin{align}\label{501}
\hat\psi^n_i(t)&:=\frac{(\theta\hat\sigma)_i\kaboom_i\sigma^2\eps}{\sum_{j=1}^I(\theta^n\hat\sigma)_i^2\kaboom_i}  V'(\theta\cdot\hat X^n(t-);\eps)\\
\label{502}
\hat\psi^n_{1,i}(t)&:=\frac{\kappa_{1,i}\sqrt{2}}{\kappa_{1,i}+\kappa_{2,i}}\hat \psi^n_i(t),\qquad
\hat\psi^n_{2,i}(t):=-\frac{\kappa_{2,i}\sqrt{2}}{(\kappa_{1,i}+\kappa_{2,i})\rho_i^{1/2}}\hat \psi^n_i(t)
,
\end{align}
and also 
\begin{align}\label{newnew3}
\psi^n_{1,i}(t):=\la_i^n+\hat\psi^n_{1,i}(t)(\la_in)^{1/2}\qquad\text{and}\qquad
\psi^n_{2,i}(t):=\mu_i^n+\hat\psi^n_{2,i}(t)(\mu_in)^{1/2}.
\end{align}
Let $\{\hat \Q^n_{j,i}\}_{j,i,n}$ be the relevant measures defined as in \eqref{216}. Notice that from Proposition \ref{prop31}, $0\le V'(\cdot;\eps)\le r$, hence it follows that all the processes mentioned in \eqref{501}--\eqref{502} are uniformly bounded by some constant. Namely, there exists a constant $C_0>0$ such that for every  $j\in\{1,2\}$, $i\in[I]$, $n\in\N$, and $t\in\R_+$,
\begin{align}
\label{new8}
|\hat\psi^n_{j,i}(t)|\le C_0.
\end{align}

We now simplify the change of measure penalty from \eqref{211}.
Since $A^n_i(\cdot)-\int_0^\cdot\hat\psi^n_{1,i}(s)ds$ is a martingale under $\hat\Q^n_{j,i}$ 
we get the following sequence of equations:
\begin{align}\label{newnew1}
\begin{split}
&L^\varrho_1(\hat \Q^n_{1,i}\|\PP^n_{1,i})\\
&\quad=\E^{\hat\Q^n_{1,i}}\Big[\int_0^\iy\rho e^{-\varrho t}\Big(\int_0^t\log\Big(\frac{\psi^n_{1,i}(s)}{\la^n_i}\Big)dA^n_i(s)-\int_0^t(\psi^n_{1,i}(s)-\la^n_i)ds\Big)dt\Big]\\
&\quad=\E^{\hat\Q^n_{1,i}}\Big[\int_0^\iy\rho e^{-\varrho t}\Big(\int_0^t\log\Big(\frac{\psi^n_{1,i}(s)}{\la^n_i}\Big)(dA^n_i(s)-\psi^n_{1,i}(s)ds)\\
&\qquad\qquad\qquad\qquad\qquad\qquad+\int_0^t\Big\{\psi^n_{1,i}(s)\log\Big(\frac{\psi^n_{1,i}(s)}{\la^n_i}\Big)-\psi^n_{1,i}(s)+\la^n_i\Big\}ds\Big)dt\Big]\\
&\quad=\E^{\hat\Q^n_{1,i}}\Big[\int_0^\iy\rho e^{-\varrho t}\Big(\int_0^t\Big\{\psi^n_{1,i}(s)\log\Big(\frac{\psi^n_{1,i}(s)}{\la^n_i}\Big)-\psi^n_{1,i}(s)+\la^n_i\Big\}ds\Big)dt\Big]\\
&\quad=\E^{\hat\Q^n_{1,i}}\Big[\int_0^\iy e^{-\varrho t}\Big\{\psi^n_{1,i}(t)\log\Big(\frac{\psi^n_{1,i}(t)}{\la^n_i}\Big)-\psi^n_{1,i}(t)+\la^n_i\Big\}dt\Big]
,
\end{split}
\end{align}
where the last equality follows by changing the order of integration. Similar calculations apply to the change of measure penalty associated with the service time. Set $y^n=\hat \psi^n_{1,i}(t)(\la_in)^{1/2}/\la^n_i$. Then,
$|y^n|\le C_1n^{-1/2}+o(n^{-1/2})$. Noticing that $|(1+y^n)\log(1+y^n)-y^n|$ is uniformly bounded over $n$ and recalling \eqref{newnew3}, we get that there exists a constant $C_2>0$, independent of $n$ and $t$, such that
\begin{align}\notag
\sum_{i=1}^I\frac{1}{\kappa_{1,i}}L^\varrho_1(\hat \Q^n_{1,i}\|\PP^n_{1,i})
+\sum_{i=1}^I\frac{1}{\kappa_{2,i}}L^\varrho_2(\hat \Q^n_{2,i}\|\PP^n_{2,i})\le C_2.
\end{align}

Fix an arbitrary sequence of controls $\{(\hat Y^n,\hat R^n)\}_n$. We focus only on those $n\in\N$ for which 
\begin{align}\label{503}
J^n(X^n(0), U^n,  R^n,\hat \Q^n;\kappa)< V(x_0;\eps)+1.
\end{align}
For the $n$'s for which the reversed inequality holds we obtain $J^n(X^n(0), U^n,  R^n,\hat \Q^n;\kappa)> V(x_0;\eps)$ and the lower bound holds without even taking a limit.
Therefore,
for every $t\in\R_+$,
\begin{align}\label{504}
e^{-\varrho t}\E^{\hat\Q^n}\left[\hat r\cdot\hat R^n(t)\right]\le \E^{\hat\Q^n}\left[\int_0^te^{-\varrho s}\hat r\cdot d\hat R^n(s)\right]< V(x_0;\eps)+C_2+1,
\end{align}
where $\hat \Q^n=\prod_{i=1}^I(\hat\Q^n_{1,i}\times\hat \Q^n_{2,i})$.
This property serves us in Lemma \ref{lem41} below when we claim tightness of a time-rescaled version of $\hat R^n$.

Recall that the intensity of $S^n_i(T^n_i)$ is $\psi^n_{2,i}U^n_i$ and that $dT^n_i(s)=U^n_i(s)ds$. Then, under the measure $\hat\Q^n$, the dynamics of $\hat X^n$ can be expressed in the following convenient form
\begin{align}\label{505}
\hat X^n_i(t)
&=\hat X^n_i(0)+\hat m^n_i t+\check A^n_i(t)-\check D^n_i(t)+\hat Y^n_i(t)-\hat R^n_i(t)\\\notag
&\qquad+\la_i^{1/2}\int_0^t\hat\psi^n_{1,i}(s)ds -\mu_i^{1/2}\int_0^{t}\hat\psi^n_{2,i}(s)dT^n_i(s),
\end{align}
where 
\begin{align}\label{new9}
\check A^n_i(t)&:=n^{-1/2}\left(A^n_i(t)-\int_0^t \psi^n_{1,i}(s)ds\right),\\
\label{new10}
\check D^n_i(t)&:=n^{-1/2}\left(S^n_i(T^n_i(t))-\int_0^{t} \psi^n_{2,i}(s)dT^n_i(s)\right).
\end{align}
Set $\check A^n=(\check A^n_1,\ldots,\check A^n_I)$ and similarly $\check D^n=(\check D^n_1,\ldots,\check D^n_I)$. By standard martingale techniques $\check A^n$ and $\check D^n$ are ${\calG}^n_t$-martingales under $\hat\Q^n$, where recall that ${\calG}^n_t$ is given in Definition \ref{def21}.(i), see e.g., the arguments given in the proof of Theorem 3.4 in \cite{kus-mar}. 

\subsubsection{Time rescaling}\label{sec421}
Recall the definition of $\theta^n$ from \eqref{thetan}. For every $n\in\N$ define
\begin{align}\label{509}
\tau^n(t):=t+\theta^n\cdot\hat R^n(t)+\theta^n\cdot\hat Y^n(t),\qquad t\in\R_+.
\end{align}
Since $\hat R^n$ and $\theta^n\cdot \hat Y^n$ are nondecreasing and RCLL (see \eqref{207}) it follows that $\tau^n$ is strictly increasing and RCLL. Moreover, for every $0\le s\le t$, 
$\tau^n(t)-\tau^n(s)\ge t-s$. The {\it time rescaled process} is given by
\begin{align}\notag
\tilde\tau^n(t):=\inf\{ s\ge 0: \tau^n(s)>t\},\qquad t\in\R_+.
\end{align}
Notice that $\tilde \tau^n$ is nondecreasing and continuous. Also, since $\hat Y^n$ is continuous and the jumps of $\hat R^n$ are of size $\frac{1}{\sqrt{n}}$, we get that there $|\tau^n(\tilde\tau^n(t))-t|\le\|\theta^n\|/\sqrt{n}$ and
\begin{align}\label{511}
0\le\tilde \tau^n(t)\le s\qquad\text{if and only if}\qquad \tau^n(s)\ge t\ge 0.
\end{align}
Define also the following rescaled processes
\begin{align}\notag
\tilde L^n(t):=\hat L^n(\tilde\tau^n(t)),
\quad \tilde A^n(t):=\check A^n(\tilde \tau^n(t)), \quad \tilde D^n(t):=\check D^n(\tilde \tau^n(t)),\quad \tilde T^n(t):=T^n(\tilde \tau^n(t)), 
\end{align}
for $L=X,Y, $ and $R$. 
\begin{lemma}\label{lem41}
%
The sequence of processes
\begin{align}\notag
\{(\check A^n,\check D^n,\tilde A^n,\tilde D^n,\tilde X^n,\tilde\tau^n,\tilde Y^n,\tilde R^n, \tilde T^n)
\}_n
\end{align}
is $\cal C$-tight. 
Let $(\check A,\check D,\tilde A,\tilde D,\tilde X,\tilde\tau,\tilde Y,\tilde R,\tilde T)$ be a limit of a weakly convergent subsequence. Then, for every $t\in\R_+$ one has,
\begin{align}\label{514}
\tilde X_i(t)&=\tilde X_i(0)+\hat m_i\tilde\tau(t)+\hat\sigma_{ii}(\theta\hat\sigma)_i\kaboom_i\int_0^tV'(\theta\cdot \tilde X(s);\eps)d\tilde\tau(s)+\hat\sigma_i\tilde B_i(t) +\tilde Y_i(t)-\tilde R_i(t),
\\\label{514b}
\tilde A(t)&=\check A(\tilde\tau(t)),\qquad \tilde D(t)=\check D(\tilde\tau(t)),\qquad  \tilde T(t):=(\rho_1,\ldots,\rho_I)\tilde\tau(t),
\end{align}
where $\tilde B=(\tilde B_1,\ldots,\tilde B_I)=\hat\sigma^{-1}(\tilde A-\tilde D)$ is a martingale w.r.t.~its own filtration and with quadratic variation 
$\tilde\tau(\cdot)\calI$, where $\calI$ is the identity matrix of order $I\times I$.
\end{lemma}
The proof is given in Section \ref{sec6}.

By reducing to a subsequence and by Skorokhod's representation theorem (see \cite[Theorem 6.7]{Bill}), we may assume without loss of generality that there is a probability space $(\Omega^* ,\calF^*,\Q^*)$ that supports 
the sequence of processes $\{(\check A^n,\check D^n,\tilde A^n,\tilde D^n,\tilde X^n,\tilde\tau^n,\tilde Y^n,\tilde R^n,\tilde T^n)\}_n$ and $\\$$(\check A,\check D,\tilde A,\tilde D,\tilde X,\tilde\tau,\tilde Y,\tilde R,\tilde T)$ such that
\begin{align}\label{521}
\limn(\check A^{n},\check D^{n},\tilde A^{n},\tilde D^{n},\tilde X^{n},\tilde\tau^{n},\tilde Y^{n},\tilde R^{n},\tilde T^{n})=(\check A,\check D,\tilde A,\tilde D,\tilde X,\tilde\tau,\tilde Y,\tilde R,\tilde T),
\end{align}
$\Q^*$-a.s., uniformly on compacts (u.o.c.)
Throughout the rest of Section \ref{sec42} we consider the probability space $(\Omega^* ,\calF^*,\Q^*)$ and w.l.o.g.~occasionally assume that \eqref{521} is at force.


The following lemma states that the time rescaled process $\tilde\tau(t)$ grows to infinity together with $t$. Loosely speaking, it implies that the processes $\{(\hat R^n,\hat Y^n)\}_n$ do not explode as $n\to\iy$ and that a right-inverse function $\tilde\tau^{-1}$ exists, see the next subsection. Hence, we can make the transformation back to the original time scale. The Lemma plays an important role in the proof of Proposition \ref{prop43} below and its proof is provided in Section \ref{sec6}.
\begin{lemma}\label{lem42}
\begin{align}\notag
\lim_{t\to\iy}\tilde\tau(t)=\iy , \quad \Q^*\text{-a.s.}
\end{align}
\end{lemma}

\subsubsection{Back to the original scale in the limiting environment}\label{sec422}
We now define the inverse of $\tilde\tau$, which brings the limit processes back to the original scale. Set
\begin{align}\notag
\tau(t):=\inf\{s\ge0 : \tilde\tau(s)>t\},\qquad t\in\R_+.
\end{align}
One can verify that $\tau$ is right-continuous and strictly increasing. Moreover, $\lim_{t\to\iy}\tau(t)=\iy$ $\Q^*$-a.s.~and from Lemma \ref{lem42}, for every $t\in\R_+$, $\tau(t)<\iy$ a.s. Finally, for every $t\in\R_+$, $\tilde \tau(\tau(t))=t$, $\tau(\tilde \tau(t))\ge t$, and
\begin{align}\notag
0\le \tilde \tau(s)\le t\qquad\text{if and only if}\qquad \tau(t)\ge s\ge 0.
\end{align}

 The time-transposed processes are defined as follows:
\begin{align}\notag
X^*(\cdot)&:=\tilde X(\tau(\cdot)),\quad A^*(\cdot):=\tilde A(\tau(\cdot)),\quad
D^*(\cdot):=\tilde D(\tau(\cdot)),\\\notag
 Y^*(\cdot)&:=\tilde Y(\tau(\cdot)),\quad
R^*(\cdot):=\tilde R(\tau(\cdot)),\quad T^*(\cdot):=\tilde T(\tau(\cdot)).
\end{align}
From \eqref{514}, the equality $\tilde \tau(\tau(t))=t$, and Lemma \ref{lem_A1}, we have for every $t\in\R_+$,
\begin{align}\label{528}
X^*_i(t)=X^*_i(0)+m_it+\hat\sigma_{ii}(\theta\hat\sigma)_i\kaboom_i\int_0^tV'(\theta\cdot X^*(s);\eps)ds+\hat\sigma_i B^*_i(t)+Y^*_i(t)-R^*_i(t),
\end{align}
where
\begin{align}\notag
B^*=( B^*_1,\ldots,B^*_I):=\hat\sigma^{-1}(A^*-D^*).
\end{align}

The relevant filtration in the limiting environment is now provided. Set
$\tilde\calG^{\sharp}_t=\tilde\calG^{\sharp}(t):=\sigma\{(\tilde X(s),\tilde A(s),\tilde D(s),\tilde Y(s),\tilde R(s)),\;0\le s\le t\}$. Notice that for every $0\le s<t<\iy$, $\{\tau(s)<t\}=\{\tilde\tau(t)>s\}\in\tilde\calG^{\sharp}(t)$. Therefore, $\tau(s)$ is an optional time for $\tilde\calG^{\sharp}(t)$. From \cite[Corollary 2.4]{kar-shr}, $\tau(s)$ is a stopping time for the complete right-continuous filtration
$\tilde \calG_t=\tilde \calG(t):=\tilde \calG^{\sharp}(t+)\vee\calN$, where $\calN$ is the collection of $\Q^*$-null sets. Using the monotonicity of $t\mapsto\tau(t)$, we get that $\calG^*_t:=\tilde\calG(\tau(t))$ is a filtration. 

\begin{proposition}\label{prop43}
The process $B^*$ is an $I$-dimensional standard Brownian motion under the filtration $\calG^*_t$.
\end{proposition}
The proof uses classical martingale arguments and is given in Section \ref{sec6}.

\subsubsection{Asymptotic lower bound}\label{sec423}
We are now ready to analyze the cost function. We start with an upper bound for the limit of the Kullback--Leibler divergences from \eqref{211}. Recall \eqref{newnew1} and consider $y^n=\hat\psi^n_{1,i}(t)(\la_in)^{1/2}/\la^n_i$. By \eqref{501}--\eqref{newnew3} and since $0\le V'\le r$ ( Proposition \ref{prop31}), $y^n\ge 0$. Notice that for every $y\ge 0$,
\begin{align}
\label{newnew23}
(1+y)\log(1+y) -y\le \frac{1}{2}y^2.
\end{align}
Since $\la_in/\la^n_i\to 1$ as $n\to\iy$, we get from the last display in \eqref{newnew1} that 
\begin{align}\notag
L^\varrho_1(\hat \Q^n_{1,i}\|\PP^n_{1,i})\le \E^{\Q^*}\Big[\frac{1}{2}\int_0^\iy e^{-\varrho t}(\hat\psi^n_{1,i}(t))^2dt\Big]+o(1)
.
\end{align}
The same calculation yields,
\begin{align}\notag
&L^\varrho_2(\hat \Q^n_{2,i}\|\PP^n_{2,i})\\\notag
&\quad\le\E^{\Q^*}\Big[\frac{1}{2}\int_0^\iy e^{-\varrho t}(\hat\psi^n_{2,i}(t))^2d\rho_i(t)\Big]+\E^{\Q^*}\Big[\int_0^\iy e^{-\vr t}(\hat\psi^n_{2,i}(t))^2d(T^n_i(t)-\rho_it)\Big]+o(1)
.
\end{align}
The last expectation above is of order $o(1)$. Indeed, by \eqref{501}--\eqref{502}, $(\hat\psi^n_{2,i}(t))^2=c_i(V'(\theta\cdot\hat X^n(t-));\eps))^2$ for some $c_i>0$. From Lemma \ref{lem_A1},
\begin{align}
\notag
\int_0^\iy e^{-\vr t}(\hat\psi^n_{2,i}(t))^2d(T^n_i(t)-\rho_it)&=\int_0^\iy e^{-\vr \tilde\tau^n(t)}(\hat\psi^n_{2,i}(\tilde\tau^n(t)))^2d(T^n_i(\tilde\tau^n(t))-\rho_i\tilde\tau^n(t))\\\notag
&=c_i\int_0^\iy e^{-\vr \tilde\tau^n(t)}(V'(\theta\cdot\tilde X^n(t-)))^2d(T^n_i(\tilde\tau^n(t))-\rho_i\tilde\tau^n(t)).
\end{align}
From Lemmas \ref{lem41} and \ref{lem_A2} the last integral converges to 0, $\Q^*$-a.s. 
Although the mentioned lemma is stated for finite time interval, the discounted cost, the boundedness of $V'$ and the bound $T^n(t)\le t$, allow us to take the integral's upper limit to be infinity. Now the uniform boundedness of the integral implies that the expectation of the last integral above converges to 0 as well. 
Using the representations above together with \eqref{300b} and \eqref{501}--\eqref{502}, it follows that
\begin{align}\label{newnew21}
&\sum_{i=1}^I\frac{1}{\kappa_{1,i}}L^\varrho_1(\hat \Q^n_{1,i}\|\PP^n_{1,i})+\sum_{i=1}^I\frac{1}{\kappa_{2,i}}L^\varrho_2(\hat \Q^n_{2,i}\|\PP^n_{2,i})\\\notag
&\qquad\le \frac{1}{2\eps}\E^{ Q^*}\Big[\int_0^\iy\eps^2\sigma^2 (V'(\theta^n\cdot\hat X^n(t);\eps))^2dt\Big]+o(1)
.
\end{align}
Notice that $\hat X^n$ has only countable number of jumps during the time interval $[0,\iy)$ and therefore we could replace $\hat X^n(t-)$ by $\hat X^n(t)$ without affecting the integral. 
From \eqref{212}, the above, and Lemma \ref{lem_A1}, one has
\begin{align}\label{549}
 &J^n( X^n(0),U^n,R^n,\hat\Q^n;\kappa)\\\notag
  &\quad\ge\E^{Q^*}\Big[\int_0^\iy e^{-\varrho t}\left\{[\hat h\cdot\hat  X^n(t)- (\eps\sigma V'(\theta^n\cdot\hat X^n(t);\eps))^2/(2\eps)]dt +\hat r\cdot  d\hat R^n(t)\right\}\Big]+o(1)
  \\\notag
  &\quad=\E^{Q^*}\Big[\int_0^\iy e^{-\varrho \tilde\tau^n(t)}\left\{[\hat h\cdot\hat  X^n(\tilde\tau^n(t))- (\eps\sigma V'(\theta^n\cdot\hat X^n(\tilde\tau^n(t));\eps))^2/(2\eps)]d\tilde\tau^n(t)\right. \\\notag&\qquad\qquad\quad\qquad\qquad\qquad\left.+\hat r\cdot  d\hat R^n(\tilde\tau^n(t))\right\}\Big]+o(1)
  \\\notag
  &\quad=\E^{Q^*}\Big[\int_0^\iy e^{-\varrho \tilde\tau^n(t)}\left\{[\hat h\cdot\tilde  X^n(t)- (\eps\sigma V'(\theta^n\cdot\tilde X^n(t);\eps))^2/(2\eps)]d\tilde\tau^n(t)  +\hat r\cdot  d\tilde R^n(t)\right\}\Big]+o(1)
  .
\end{align}
Recall that $\hat h\cdot\tilde X^n$ and $\tilde R^n$ are nonnegative, that $\tilde X^n$ is uniformly bounded, and also that $\tilde R^n$ is nondecreasing, continuous, and bounded on any compact time interval. Then from \eqref{521}
,  Lemma \ref{lem_A2}, and the bounded convergence theorem, we get that $\Q^*$-almost-surely for every $s\in\R_+$, one has 
\begin{align}\notag
&\liminf_{n\to\iy}\;
\int_0^\iy  e^{-\varrho \tilde\tau^n(t)}\left\{\hat h\cdot\tilde  X^n(t)d\tilde\tau^n(t)  +\hat r\cdot  d\tilde R^n(t)\right\}\\\notag
&\quad\ge \limn\;
\int_0^s  e^{-\varrho \tilde\tau^n(t)}\left\{\hat h\cdot\tilde  X^n(t)d\tilde\tau^n(t)  +\hat r\cdot  d\tilde R^n(t)\right\}\\\notag
&\quad=\int_0^s  e^{-\varrho \tilde\tau(t)}\left\{\hat h\cdot\tilde  X(t)d\tilde\tau(t)  +\hat r\cdot  d\tilde R(t)\right\}.
\end{align}
Taking $s\to\iy$ (the upper limit of the integral) first and then $\E^{\Q^*}$ on both sides, we obtain
\begin{align}\notag
&\liminf_{n\to\iy}\;\E^{\Q^*}\Big[
\int_0^\iy  e^{-\varrho \tilde\tau^n(t)}\left\{\hat h\cdot\tilde  X^n(t)d\tilde\tau^n(t)  +\hat r\cdot  d\tilde R^n(t)\right\}\Big]\\\notag
&\quad\ge\E^{\Q^*}\Big[\int_0^\iy  e^{-\varrho \tilde\tau(t)}\left\{\hat h\cdot\tilde  X(t)d\tilde\tau(t)  +\hat r\cdot  d\tilde R(t)\right\}\Big].
\end{align}
Using now $0\le V'\le r$, $(\tilde X^n,\tilde \tau^n)\To(\tilde X,\tilde\tau)$,  Lemma \ref{lem_A2}, and the bounded convergence theorem, we similarly get that
\begin{align}\notag
&\liminf_{n\to\iy}\;\E^{\Q^*}\Big[
-\int_0^\iy  e^{-\varrho \tilde\tau^n(t)} (\eps\sigma V'(\theta^n\cdot\tilde X^n(t);\eps))^2/(2\eps)d\tilde\tau^n(t)\Big]\\\notag
&\quad\ge\E^{\Q^*}\Big[-\int_0^\iy  e^{-\varrho \tilde\tau(t)} (\eps\sigma V'(\theta^n\cdot\tilde X(t);\eps))^2/(2\eps)d\tilde\tau(t) \Big].
\end{align}
Combining the last two bounds, 
\begin{align}\label{552}
&\underset{n\to\iy}{\liminf}\;\E^{\Q^*}\Big[\int_0^\iy  e^{-\varrho \tilde\tau^n(t)}\left\{[\hat h\cdot\tilde  X^n(t)- (\eps\sigma V'(\theta^n\cdot\tilde X^n(t);\eps))^2/(2\eps)]d\tilde\tau^n(t)  +\hat r\cdot  d\tilde R^n(t)\right\}\Big]
\\\notag
&\qquad\ge \E^{\Q^*}\Big[\int_0^\iy e^{-\varrho \tilde\tau(t)}\left\{[\hat h\cdot\tilde  X(t)- (\eps\sigma V'(\theta\cdot\tilde X(t);\eps))^2/(2\eps)]d\tilde\tau(t)  +\hat r\cdot  d\tilde R(t)\right\}\Big]\\\notag
&\qquad= \E^{\Q^*}\Big[\int_0^\iy e^{-\varrho t}\left\{[\hat h\cdot  X^*(t)- (\eps\sigma V'(\theta\cdot X^*(t);\eps))^2/(2\eps)]dt  +\hat r\cdot  d R^*(t)\right\}\Big]\\\notag
&\qquad\ge \E^{\Q^*}\Big[\int_0^\iy e^{-\varrho t}\left\{[ h(\theta\cdot  X^*(t))- (\eps\sigma V'(\theta\cdot X^*(t);\eps))^2/(2\eps)]dt  + r\cdot  d(\theta\cdot R^*(t))\right\}\Big],
\end{align}
where we used Lemma \ref{lem_A1} to get the equality. Indeed, recall that $\tilde\tau(\tau(t))=t$, $X^*(t)=\tilde X(\tau(t))$, $R^*(t)=\tilde R(\tau(t))$, and Lemma \ref{lem42}. The last inequality follows since by the definitions of $h$ and $r$ (see \eqref{319a}--\eqref{319b} and \cite[(2.45)--(2.46)]{Cohen2017}),
$
h(\theta\cdot X^*(t))\le\hat h\cdot X^*(t)$ and $
\int_0^\iy e^{-\varrho t} r d( \theta\cdot R^*(t))\le \int_0^\iy e^{-\varrho t}\hat r\cdot dR^*(t)$.

Denote $B:=(\sum_{i=1}^I(\theta\hat\sigma)_iB^*_i)/\sigma$ and $(X,Y,R):=(\theta\cdot X^*,\theta\cdot Y^*,\theta\cdot R^*)$. Then, from  \eqref{514}, together with \eqref{300b}, \eqref{310b}, and the limit $x_0=\limn \theta^n\cdot\hat X^n(0)$, we have that
\begin{align}\label{asaf002}
X(t)=x_0+mt+\sigma\int_0^t\eps\sigma V'(X(s);\eps)ds+\sigma B(t)+Y(t)-R(t),\quad t\in\R_+.
\end{align}
From Proposition \ref{prop43} we get that $B$ is a standard one-dimensional Brownian motion w.r.t.~$\calG^*_t$. Hence, Property (i) of Definition \ref{def32} holds. Also, Property (ii) follows since $\theta^n\cdot \hat Y^n$ and $\hat R^n$ are nonnegative and nondecreasing, See \eqref{207}. Hence, \eqref{321} and the definition of $\Q_V$ in Proposition \ref{prop31} imply that the last expectation in the sequence of relations in \eqref{552} equals $J(x_0,Y,R,\Q_V;\eps)$. Together with \eqref{549}  we have,
\begin{align}\notag
\underset{n\to\iy}{\liminf}\;J^n( X^n(0),U^n,R^n,\hat\Q^n;\kappa)\ge J(x_0,Y,R,\Q_V;\eps)\ge V(x_0;\eps).
\end{align}
Since the sequence of policies $\{(U^n,R^n)\}_n$ is arbitrary, it follows that
\begin{align}\notag
\underset{n\to\iy}{\liminf}\;V^n( X^n(0);\kappa)\ge V(x_0;\eps).\end{align}

\hfill$\Box$

\subsection{Proof of \eqref{new5}}\label{sec43}

In this section we consider a sequence of arbitrary strategies for the maximizer and show that in the limit, the candidate policy for the DM is bounded above by the value function of the RSDG. That is, asymptotically, the maximizer cannot punish the DM more than the equilibrium in the RSDG. 

Recall that the maximizer's strategy can be equivalently formulated as rate perturbations. We start with considering an arbitrary sequence of strategies for the maximizer that is not too costly. That is, the rates will be only moderately perturbed by the maximizer. Then in Section \ref{sec432} we prove a state-space collapse property. Specifically, we show that by using the candidate policy, the dynamics of the buffers' sizes stay close to $\gamma^a$ from \eqref{asaf003} under the perturbed rates. In Section \ref{sec433} we asymptotically bound the expected cost; in order to estimate the change of measure penalty, we truncate the rate processes $\{\psi^n_{j,i}\}_{n,j,i}$ and show that by doing this the penalty does not change much.

%
%
%
%
%
%
%
%
%
%
\subsubsection{The maximizer's perspective}\label{sec431}
Consider an arbitrary sequence of measures chosen by the maximizer in the QCP, $\{\hat \Q^n\}_{n\in\N}$, where each $\hat \Q^n\in\hat \calQ^n( X^n(0))$.
Recall that every measure $\hat \Q^n\in\hat \calQ^n( X^n(0))$ is associated with the processes $\{\psi^n_{j,i}\}_{j,i}$, see \eqref{216}. These processes stand for the `new' intensity of the processes $\{A^n_i\}_{i,n}$ and $\{S^n_i(T^n_i)\}_{i,n}$. To simplify the notation and some of the arguments, we consider one probability space $(\Omega,{\calG},\hat\Q)$ that supports
the processes $\{(A^n_i,S^n_i(T^n_i))\}_{i,n}$ and under-which, for every $n\in\N$, the relevant intensities are $\{\psi^n_{j,i}\}_{j,i,n}$. However, occasionally, when we want to emphasize the relevant measure, we use the measures $\{\hat\Q^n_{j,i}\}_{n,j,i}$ and $\{\hat\Q^n\}_n$.

Notice that the same calculation given in \eqref{newnew1} is valid here as well, so
\begin{align}\label{newnew2}
\begin{split}
L^\varrho_1(\hat \Q^n_{1,i}\|\PP^n_{1,i})
&=\E^{\hat\Q^n_{1,i}}\Big[\int_0^\iy e^{-\varrho t}\Big(\psi^n_{1,i}(t)\log\Big(\frac{\psi^n_{1,i}(t)}{\la^n_i}\Big)-\psi^n_{1,i}(t)+\la^n_i\Big)dt\Big],\\
L^\varrho_2(\hat \Q^n_{2,i}\|\PP^n_{2,i})&=\E^{\hat\Q^n_{2,i}}\Big[\int_0^\iy e^{-\varrho t}\Big(\psi^n_{2,i}(t)\log\Big(\frac{\psi^n_{2,i}(t)}{\mu^n_i}\Big)-\psi^n_{2,i}(t)+\mu^n_i\Big)dT^n_i(t)\Big].
\end{split}
\end{align}

We now show that without any loss, the maximizer can be restricted to measures $\hat \Q^n\in\hat \calQ^n( X^n(0))$, which are `not too far away' from the reference measure $\PP^n$. The idea behind it, as will be provided rigorously in the proof below, is that by changing the rate, the rejection cost will contribute at most a linear cost, while the penalty for the change of measure is super-linear. Similarly to the bound in \eqref{503}, without loss of generality, we may and will assume that
\begin{align}\label{new13}
J^n( X^n(0),U^n(a),R^n(a),\hat\Q^n;\kappa)\ge V(x_0;\eps)-1.
\end{align}

Define 
\begin{align}
\label{newnew22}
\hat\psi^n_{1,i}(t):=(\la_i n)^{-1/2}\left(\psi^n_{1,i}(t)-\la^n_i\right),\qquad\text{and}\qquad\hat\psi^n_{2,i}(t):=(\mu_i n)^{-1/2}\left(\psi^n_{2,i}(t)-\mu^n_i\right).
\end{align}
Recall that in the previous subsection the maximizer was given a specific strategy for which $\sup_{t,j,i}|\hat\psi^n_{j,i}(t)|$ was bounded. Since we  consider now a sequence of arbitrary strategies for the maximizer, it does not hold in this case. However, as we show in Proposition \ref{lem44} and in Section \ref{sec433}, this is approximately the case.
\begin{proposition}\label{lem44}
There exists $M>0$ 
such that for every $n\ge n_0$ and every 
$i\in[I]$, 
\begin{align}\notag
&\E^{\hat\Q}\Big[\int_0^\iy e^{-\varrho t}\Big(
 |\hat\psi^n_{1,i}(t)|dt+ |\hat\psi^n_{2,i}(t)|dT^n_i(t)\Big)\\\notag
&\qquad+\int_0^\iy e^{-\varrho t}\Big(
(\hat\psi^n_{1,i}(t))^2dt+ (\hat\psi^n_{2,i}(t))^2dT^n_i(t)\Big)\Big]\le M
\end{align}
and
\begin{align}\notag
 &\E^{\hat \Q}\Big[
 \int_0^\iy e^{-\varrho t}\Big(\psi^n_{1,i}(t)\log\Big(\frac{\psi^n_{1,i}(t)}{\la^n_i}\Big)-\psi^n_{1,i}(t)+\la^n_i\Big)dt\\\notag
 &\qquad+\int_0^\iy e^{-\varrho t}\Big(\psi^n_{2,i}(t)\log\Big(\frac{\psi^n_{2,i}(t)}{\mu^n_i}\Big)-\psi^n_{2,i}(t)+\mu^n_i\Big)dT^n_i(t)\Big]
 \le M.
\end{align}
\end{proposition}
{\bf Proof.} The second bound follows from the first one together with inequality \eqref{newnew23}. Thus, we only prove the first one. 
As argued in \cite{ata-shi}, at the bottom of page 595, for every $t\in\R_+$,
\begin{align}\notag
\|\hat R^n(t)\|\le C(1+t+\|\hat A^n\|_t+\|\hat D^n\|_t),
\end{align}
where in the above expression, and in the rest of the proof, $C$ refers to a finite positive constant
that is independent of $n$ and $t$ and which can change from one line to the next. Denote
\begin{align}\notag
\varphi^n(t):=\sum_{i=1}^I\left(|\hat\psi^n_{1,i}(t)|
+|\hat\psi^n_{2,i}(t)|U^n_i(t)\right), \qquad t\in\R_+.
\end{align}
Recall the definitions of $\check A^n, \check D^n$ given in \eqref{new9} and \eqref{new10}, then
\begin{align}\label{new19}
\|\hat R^n(t)\|\le C\Big(1+t+\int_0^t\varphi^n(s)ds+\|\check A^n\|_t+\|\check  D^n\|_t\Big).
\end{align}
Applying Burkholder--Davis--Gundy inequality to $\check A^n$ and $\check D^n$, and noticing that 
$n^{-1}|\psi^n_{j,i}|\le C(1+n^{-1/2}|\hat\psi^n_{j,i}|)$,
we have
\begin{align}
\label{newnew14}
\begin{split}&\E^{\hat \Q}[\|\check A^n\|_t^2]\le Cn^{-1}\E^{\hat \Q}\Big[\int_0^t\psi^n_{1,i}(s)ds\Big]\le C\Big(t+n^{-1/2}\E^{\hat \Q}\Big[\int_0^t\hat \psi^n_{1,i}(s)ds\Big]\Big),\\
&\E^{\hat \Q}[\|\check D^n\|_t^2]\le Cn^{-1}\E^{\hat \Q}\Big[\int_0^t\psi^n_{2,i}(s)dT^n_i(s)\Big]\le C\Big(t+n^{-1/2}\E^{\hat \Q}\Big[\int_0^t\hat \psi^n_{2,i}(s)dT^n_i(s)\Big]\Big).
\end{split}
\end{align}
Hence,
\begin{align}\label{new12}
\E^{\hat\Q}\left[\|\hat R^n(t)\|\right]\le C\Big(1+t+\E^{\hat\Q}\Big[\int_0^t\varphi^n(s)ds\Big]\Big).
\end{align}
An application of integration by parts yields
\begin{align}\label{new11}
\int_0^\iy e^{-\varrho t}\hat r\cdot d\hat R^n(t)&=\left[-\rho^{-1}e^{-\varrho t}\hat r\cdot \hat R^n(t)\right]_{t=0}^{t=\iy}+\int_0^\iy \varrho e^{-\varrho t}\hat r\cdot \hat R^n(t)dt,
\end{align}
and changing the order of integration implies
\begin{align}\notag
\int_0^\iy e^{-\varrho t}\int_0^t\varphi^n(s)ds
&=\int_0^\iy \varrho e^{-\varrho t}\varphi^n(t)dt.
\end{align}
Notice that the first term in \eqref{new11} is non-positive. Taking expectation on both sides of it and using the bound \eqref{new12}, we get that
\begin{align}\notag
\E^{\hat\Q}\Big[\int_0^\iy e^{-\varrho t}\hat r\cdot d\hat R^n(t)\Big]
&\le C\Big(1+\E^{\hat\Q}\Big[\int_0^\iy \varrho e^{-\varrho t}\Big\{\int_0^t\varphi^n(s)ds\Big\}dt\Big]\Big)\\\notag
&= C\Big(1+\E^{\hat\Q}\Big[\int_0^\iy \varrho e^{-\varrho t}\varphi^n(t)dt\Big]\Big).
\end{align}
Clearly, $\hat h\cdot\hat X^n$ is bounded. Then \eqref{new13} and the last bound yield that
\begin{align}\label{new14}
\begin{split}
&\E^{\hat \Q}\Big[\int_0^\iy e^{-\varrho t}\left((\hat\psi^n_{1,i}(t))^2dt
+(\hat\psi^n_{2,i}(t))^2dT^n_i(t)\right)\Big]\\
&\quad
\le C_0\left(1+\E^{\hat \Q}\Big[\int_0^\iy e^{-\varrho t}\left(|\hat\psi^n_{1,i}(t)|dt+|\hat\psi^n_{2,i}(t)|dT^n_i(t)\right)\Big]\right),
\end{split}
\end{align}
for some $C_0>0$, independent of $n$.  
%
%
%
Since the function $z\mapsto z^2$ is super-linear, there is a constant $C_1\in\R$ such that for any $z\in\R$, $z^2\ge C_1+2C_0z$, and thus 
\begin{align}\notag
\begin{split}&\E^{\hat \Q}\Big[\int_0^\iy e^{-\varrho t}\left((\hat\psi^n_{1,i}(t))^2dt+(\hat\psi^n_{2,i}(t))^2dT^n_i(t)\right)\Big]\\
&\quad
\ge C_1+2C_0\left(\E^{\hat \Q}\Big[\int_0^\iy e^{-\varrho t}\left(|\hat\psi^n_{1,i}(t)|dt+|\hat\psi^n_{2,i}(t)|dT^n_i(t)\right)\Big]\right).
\end{split}
\end{align}
Together with \eqref{new14}, we obtain that,
\begin{align}\notag
\E^{\hat \Q}\Big[\int_0^\iy e^{-\varrho t}\Big(
|\hat\psi^n_{1,i}(t)|dt+ |\hat\psi^n_{2,i}(t)|dT^n_i(t)\Big)\Big]\le 1-C_1/C_0,
\end{align}
and the result holds by another application of \eqref{new14} and the bound above. 

\qed

Define 
\begin{align}
\notag
\hat\Psi^n_{1,i}(\cdot):=\int_0^\cdot\hat\psi^n_{1,i}(t)dt,\qquad\text{and}\qquad\hat\Psi^n_{2,i}(\cdot):=\int_0^\cdot\hat\psi^n_{2,i}(t)dT^n_i(t),
\end{align}
and set $\hat\Psi^n_j=(\hat\Psi^n_{j,i};i\in[I])$. 
The next lemma assures that $\{(\hat\Psi^n_1,\hat\Psi^n_2)\}_n$ has a converging subsequence and that any of its limit points has continuous paths with probability one.
\begin{lemma}\label{lem45b}
The sequence $\{(\hat\Psi^n_1,\hat\Psi^n_2)\}_n$ is $\calC$-tight.
\end{lemma}
The proof is given in Section \ref{sec6}.

%

\subsubsection{State-space collapse - staying close to the minimizing curve}\label{sec432}
We consider a set of one-dimensional processes. It is obtained by multiplying the scaled processes by $\theta^n$ (see \eqref{thetan}). 
For its definition denote 
\begin{align}
\label{newnew16}
\hat W^n:=\check A^n-\check D^n+\hat m^n\cdot,\end{align} and 
\begin{align}
\label{newnew16b}
  W^{\sharp,n}:=\theta^n\cdot(\hat W^n),\quad
  X^{\sharp,n}:=\theta^n\cdot\hat X^n,\quad
  Y^{\sharp,n}:=\theta^n\cdot\hat Y^n,\quad
  R^{\sharp,n}:=\theta^n\cdot\hat R^n.
\end{align}
Moreover, set $ \psi^{\sharp,n}_1:=\sigma^{-1}\sum_{i=1}^I  \theta^n_i\la_i^{1/2}\hat\psi^n_{1,i}$, $\;  \psi^{\sharp,n}_2:=\sigma^{-1}\sum_{i=1}^I \theta^n_i\mu_i^{1/2}\hat\psi^n_{2,i}$, and
 \begin{align}\label{new42}
 \Psi^{\sharp,n}_1(\cdot):=\int_0^\cdot \psi^{\sharp,n}_1(s)ds, \qquad \Psi^{\sharp,n}_2(\cdot):=\int_0^\cdot \psi^{\sharp,n}_2(s)dT^n_i(s).
 \end{align}
The identity from \eqref{505} is valid here as well and can be expressed as 
\begin{align}\label{newnew24}
\hat X^n_i(t)
&=\hat X^n_i(0)+ \hat m^n_i t+\check A^n_i(t)-\check D^n_i(t)+\hat Y^n_i(t)-\hat R^n_i(t)+\la_i^{1/2}\hat\Psi^n_{1,i}(t) -\mu_i^{1/2}\hat\Psi^n_{2,i}(t).
\end{align}
As a result,
\begin{align}
  \label{34}
  X^{\sharp,n}(t)=X^{\sharp,n}(0)+W^{\sharp,n}(t)+Y^{\sharp,n}(t)-R^{\sharp,n}(t)+\sigma(\Psi^{\sharp,n}_1(t)-\Psi^{\sharp,n}_2(t)),\quad t\in\R_+.
\end{align}

We now state a couple of results from \cite{ata-shi} that are needed in the next subsection.  The proofs in our case are almost identical, yet require some technical modifications and they are deferred to Section \ref{sec6}. The next two lemmas are needed to prove Proposition \ref{lem47} and Lemma \ref{lem48} below.
The next lemma states that under the candidate policy given in Section \ref{sec41}, the busy time of every buffer converges to its traffic intensity. 
\begin{lemma}\label{lem45}
For every $i\in[I]$, \{$T^n_i\}_n$ converges u.o.c.~to $\bar T_i$, where $\bar T_i(t):=\rho_it$, $t \in\R_+$. Moreover, $\{(\check A^n,\check D^n,\hat W^n)\}_n$ is $\calC$-tight.
\end{lemma}

We now define another set of processes generated from the last set. Let $\tau^n$ be the first time a forced rejection occurred in the $n$-th system and set  
\begin{align}
\begin{split}
\label{newnew16c}
W^{\circ,n}(\cdot)&:=W^{\sharp,n}(\cdot\w\tau^n),\qquad X^{\circ,n}(\cdot):=X^{\sharp,n}(\cdot\w\tau^n),\qquad
Y^{\circ,n}(\cdot):=Y^{\sharp,n}(\cdot\w\tau^n),\\
 R^{\circ,n}(\cdot)&:=R^{\sharp,n}(\cdot\w\tau^n),\qquad
\Psi_1^{\circ,n}(\cdot):=\Psi_1^{\sharp,n}(\cdot\w\tau^n),\qquad \Psi_2^{\circ,n}(\cdot):=\Psi_2^{\sharp,n}(\cdot\w\tau^n).
\end{split}
\end{align} 
A minor modification from \cite{ata-shi} that can be observed immediately 
is that
we defined $\hat W^n$ and $W^{\sharp,n}$
using $\check A^n$ and $\check D^n$ instead of $\hat A^n$ and $\hat D^n$. 

We continue the proof of \eqref{new5} with the case that the initial state lies close to the minimizng curve. That is,
\begin{align}\label{new23}
\lim_{n\to\iy} \hat X^n(0)-\gamma^a(x_0)=0 \qquad\text{ and }\qquad \theta^n\cdot\hat X^n(0)\in[0,\aey ],\quad n\in\N,
\end{align}
where recall that $x_0:=\theta^n\cdot\hat X^n(0)$.
As argued in \cite[Step 5, p.~596]{ata-shi}, this condition can be relaxed by considering a stopping time that indicates when the state is `sufficiently close' to the minimizing path. This stopping time converges to zero. Then we may continue from that stopping time, in the same way. The proof follows by the same lines as in case that  
the initial state lies close to the minimizing curve but with heavier notation and therefore omitted.

A modification of the two next results appear in \cite[Step 1, pp.~586--588]{ata-shi}.

\begin{lemma}\label{lem46}
For every $n\in\N$,
\begin{align}\label{new32}
(\aey \wedge X^{\circ,n},Y^{\circ,n},R^{\circ,n})=\Gamma_{[0,\aey ]}(X^{\sharp,n}(0)+W^{\circ,n}+\sigma(\Psi^{\circ,n}_1-\Psi^{\circ,n}_2)+E^n),
\end{align}
where
\begin{align}\label{new34}
E^n:=(\aey \wedge X^{\circ,n})-X^{\circ,n}\To 0.
\end{align}
Moreover, the sequence $\{(W^{\circ,n},X^{\circ,n},Y^{\circ,n},R^{\circ,n},\Psi^{\circ,n}_1,\Psi^{\circ,n}_2)\}_n$ is $\calC$-tight.
\end{lemma}

Next we provide a state space collapse result and claim that the multidimensional process $\hat X^n$ lies close to the minimizing curve. Recall the dependency of the candidate policy on the parameter $\delta_0$. 
\begin{proposition}\label{lem47}
For sufficiently small $\delta_0$, the following limit holds u.o.c.
\begin{align}\notag
 \hat X^n-\gamma^a(X^{\sharp,n})\overset{\hat\Q}{\Longrightarrow}0,
\end{align}
as $n\to\iy$ and moreover for every $T>0$, $\limn\hat\Q(\tau^n<T)=0$.
\end{proposition}

From the limit $\lim_n\hat\Q(\tau^n<T)\to0$, we conclude that $\{(W^{\sharp,n},X^{\sharp,n},Y^{\sharp,n},R^{\sharp,n},\Psi^{\sharp,n}_1,\Psi^{\sharp,n}_2)\}_n$ is $\calC$-tight.

\subsubsection{Asymptotic bound for the costs}\label{sec433}

In this subsection we take advantage of the convergence of the dynamics to the minimizing curve we just argued. We start by showing that the expected holding and rejection components of the cost can be approximated by equivalent components  associated with one dimensional dynamics. Then we approximate the one-dimensional dynamics with simpler dynamics for which the intensities used by the maximizer are truncated in some sense. The difference between the expected holding and rejection cost components of these dynamics are shown to be small. The reason for this reduction is as follows. Although we know that the $\hat\Q$-a.s.~absolutely continuous valued processes $(\hat\Psi^n_1,\hat\Psi^n_2)$ are $\calC$-tight and therefore have a converging subsequence whose limit is denoted by $(\hat\Psi_1,\hat\Psi_2)$, it does not imply that $\hat\Q$-a.s.~the paths of $(\hat\Psi_1,\hat\Psi_2)$ are absolutely continuous. Hence, we cannot argue for example that $\hat\Psi_1$ is of the form $\hat\Psi_{1,i}=\int_0^\cdot\hat\psi_{1,i}(t)dt$ for some $\hat\psi_{1,i}$. As a consequence, we cannot express the limiting measure for the maximizer, nor the change of measure penalty using $\hat\psi_{1,i}$ as can  be done for example in \eqref{newnew2} by substituting $\psi^n_{1,i}=\la^n_i+\hat\psi^n_{1,i}(\la_in)^{1/2}$. After this reduction we bound the change of measure penalty. 
Finally, the expected cost associated with the uniformly bounded rate dynamics  is approximated by the value function of the RSDG.

\skp\noi{\bf One dimensional reduction.} 
We start with showing the following uniform bound
\begin{align}\label{new26}
\limsup_{n\to\iy} \E^{\hat\Q}\int_0^\iy e^{-\varrho t}(\hat R^n(t))^2dt<\iy.
\end{align}
To establish this, recall the bounds in \eqref{new19}--\eqref{newnew14}. As a result, 
\begin{align}\notag
\E^{\hat\Q}\left[\|\hat R^n(t)\|^2\right]\le C\Big(1+t^2+\E^{\hat\Q}\Big[\|\hat\Psi^n_1(t)\|^2+\|\hat\Psi^n_2(t)\|^2\Big]\Big),
\end{align}
for some constant $C>0$ independent of $n$ and $t$. Together with Proposition \ref{lem44},  \eqref{new26} follows.
%
%

Now turn to the holding and rejection costs.
By the definitions of the rejection mechanism and of $R^{\sharp,n}$ in Section \ref{sec432}, on the event $\{\tau^n>T\}$, one has
\begin{align}\notag
\hat R^n(T)=\hat R^n_{i^*}(T)e_{i^*}=\mu_{i^*} R^{\sharp,n}(T)e_{i^*}.
\end{align}

Relations \eqref{new12} and \eqref{new11} imply that 
\begin{align}\label{newnew28}
\E^{\hat\Q}\Big[\int_0^\iy e^{-\vr t}\hat r\cdot d\hat R^n(t)\Big]=\E^{\hat\Q}\Big[\int_0^\iy \vr e^{-\vr t}\hat r\cdot\hat R^n(t)dt\Big].
\end{align}
Recall that $h^a(w)=\hat h\cdot\gamma^a(w)$. Then, Proposition \ref{lem47}, the boundedness of $\hat X^n$ and $h^a(X^{\sharp,n})$, the bound $ R^{\sharp,n}\le\|\theta^n\|\|\hat R^n\|$, and \eqref{new26} imply that
\begin{align}\notag
\limn\;\E^{\hat\Q}\Big[\int_0^\iy e^{-\varrho t}\{\hat h\cdot \hat X^n(t)+\varrho\hat r\cdot\hat R^n(t)\}dt
-
\int_0^\iy e^{-\varrho t}\{h^a(X^{\sharp,n}(t))+\varrho r\hat R^{\sharp,n}(t)\}dt\Big]=0,
\end{align}
where we used the identity $r=r_{i^*}\mu_{i^*}$, see \eqref{320}.
Using the limit $\limn\hat\Q(\tau^n<T)\to0$ and \eqref{newnew28}, we deduce,
\begin{align}\label{new40}
\begin{split}
\limn\;\E^{\hat\Q}\Big[&\int_0^\iy e^{-\varrho t}\{\hat h\cdot \hat X^n(t)dt+\hat r\cdot d\hat R^n(t)\}\\
&-
\int_0^\iy e^{-\varrho t}\{h^a(X^{\circ,n}(t))+\varrho r\hat R^{\circ,n}(t)\}dt\Big]=0.
\end{split}
\end{align}

\skp\noi{\bf Truncated intensities.} 
We now show that the maximizer can use probability measures for which $\{\hat\psi^n_{j,i}\}_{j,i,n}$ are uniformly bounded from above by a sufficiently large constant $k$ without too much lost.

Recall that under $\hat\Q$,
$A^n_i(\cdot)-\int_0^\cdot\psi^n_{1,i}(t)dt$  and $S^n_i(T^n_i(\cdot))-\int_0^\cdot\psi^n_{2,i}(t)dT^n_i(t)$ are martingales. We now truncate the intensities and consider the $\hat\Q$ martingales
 $A^{n,k}_i(\cdot)-\int_0^\cdot\psi^{n,k}_{1,i}(t)dt$  and $S^{n,k}_i(T^{n,k}_i(\cdot))-\int_0^\cdot\psi^{n,k}_{2,i}(t)dT^{n,k}_i(t)$, where
 \begin{align}\notag
\psi^{n,k}_{1,i}(t)&:=\psi^n_{1,i}(t)-(\la_in)^{1/2}\hat\psi^n_{1,i}(t)\one_{\{|\hat\psi^n_{1,i}(t)|>k\}},\\\notag
\psi^{n,k}_{2,i}(t)&:=\psi^n_{2,i}(t)-(\mu_in)^{1/2}\hat\psi^n_{2,i}(t)\one_{\{|\hat\psi^n_{2,i}(t)|>k\}},
\end{align}
and $T^{n,k}=(T^{n,k}_i:i\in[I])$ is the DM's control associated with the intensities $\{\psi^{n,k}_{j,i}\}_{j,i}$.
Denote $\hat\psi^{n,k}_{j,i}(\cdot):=\hat\psi^n_{j,i}(\cdot)\one_{\{|\hat\psi^n_{j,i}(\cdot)|\le k\}}$, $j\in\{1,2\}$. Clearly, $|\hat\psi^{n,k}_{j,i}|\le k$ and
 \begin{align}\notag
\psi^{n,k}_{1,i}(t)=\la_i^n+\hat\psi^{n,k}_{1,i}(t)(\la_in)^{1/2}+o(n^{1/2}),\\\notag
\psi^{n,k}_{2,i}(t)=\mu_i^n+\hat\psi^{n,k}_{2,i}(t)(\mu_in)^{1/2}+o(n^{1/2}),
\end{align}
The processes $ (A^{n,k},D^{n,k})$ and $ (A^n,D^n)$ are coupled such that for every $i\in[I]$,
\begin{align}\label{new36}
\begin{split}
&(A^n_i-A^{n,k}_i)(\cdot)- (\la_in)^{1/2}\int_0^\cdot\hat\psi^{n}_{1,i}(t)\one_{\{|\hat\psi^n_{1,i}(t)|>k\}}dt,\\
&(D^n_i-D^{n,k})_i(\cdot)- (\mu_in)^{1/2}\int_0^\cdot\hat\psi^{n}_{2,i}(t)\one_{\{|\hat\psi^n_{2,i}(t)|>k\}}dT^{n,k}_i(t)
\end{split}
\end{align}
are martingales.

Define $\check A^{n,k}:=(\check A^{n,k}_i:i\in[I])$, with $\check A^{n,k}_i:=n^{-1/2}(A^{n,k}_i(\cdot)-\int_0^\cdot\psi^{n,k}_{1,i}(t)dt)$ and similarly define $\check D^{n,k}$. 
For every $L\in\{W,X,Y,R,\Psi_1,\Psi_2\}$ let $\hat L^{n,k}$ and $L^{\sharp,n,k}$ be defined as $\hat L^n$ and $L^{\sharp,n}$, where the intensities $(\psi^{n,k}_{1,i},\psi^{n,k}_{2,i})$ are replacing $(\psi^n_{i,i},\psi^n_{2,i})$, $i\in[I]$ and let $T^{n,k}_i$ be the equivalence of $T^n_i$ in this setup. Also, let $\tau^{n,k}$ be the first time a forced rejection occurs in this setup and set $L^{\circ,n,k}:=L^{\sharp,n,k}(\cdot\wedge\tau^{n,k})$.
Clearly, Lemmas \ref{lem45}, \ref{lem46}, and Proposition \ref{lem47} hold in this case as well, where the superscript $n$ is replaced by $n,k$. For the sake of exposition, we state all the necessary results here.

As a private case of Lemma \ref{lem45} we get that $\{T^{n,k}\}_n$ converges u.o.c.~to $\bar T=(\bar T_1,\ldots,\bar T_I)$, where recall that $\bar T_i(t)=\rho_i t$, $t\in\R_+$.
Proposition \ref{lem47} implies that for sufficiently small $\delta_0$, the following limits hold
\begin{align}\notag
\limn \hat X^{n,k}-\gamma^a(X^{\sharp,n,k})\overset{\hat\Q}{\Longrightarrow}0,
\end{align}
and for every $k,T>0$,
\begin{align}\label{new43}
\hat\Q(\tau^{n,k}<T)=0.
\end{align}
 \begin{lemma}\label{lem48}
For every $n\in\N$,
\begin{align}\label{new33}
(\aey \wedge X^{\circ,n,k},Y^{\circ,n,k},R^{\circ,n,k})=\Gamma_{[0,\aey ]}(X^{\sharp,n}(0)+W^{\circ,n,k}+\sigma(\Psi^{\circ,n,k}_1-\Psi^{\circ,n,k}_2)+E^{n,k}),
\end{align}
where
\begin{align}\label{new35}
E^{n,k}:=(\aey \wedge X^{\circ,n,k})-X^{\circ,n,k}\To 0.
\end{align}
Moreover, the sequence $\{(W^{\circ,n,k},X^{\circ,n,k},Y^{\circ,n,k},R^{\circ,n,k},\Psi^{\circ,n,k}_1,\Psi^{\circ,n,k}_2,\hat\Psi^{n,k}_1,\hat\Psi^{n,k}_2)\}_n$ is $\calC$-tight and any subsequential limit of it $(\bar W^{\circ,k},\bar X^{\circ,k},\bar Y^{\circ,k},\bar R^{\circ,k},\bar \Psi^{\circ,k}_1,\bar \Psi^{\circ,k}_2,\hat\Phi^k_1,\hat\Phi^k_2)$ satisfies $\hat\Q$-a.s.
\begin{align}\label{new35a}
(\bar X^{\circ,k},\bar Y^{\circ,k},\bar R^{\circ,k})=\Gamma_{[0,\aey ]}(\bar X^{\circ,k}(0)+\bar W^{\circ,k}+\sigma(\bar \Psi^{\circ,k}_1-\bar \Psi^{\circ,k}_2)),
\end{align}
where $\bar X^{\circ,k}(0)= X^{\sharp,n}(0)=x_0$,  $\bar W^{\circ,k}$ is an $(m,\sigma)$-Brownian motion w.r.t.~the filtration $\calF_t:=\sigma\{\bar W^{\circ,k}(s),\bar X^{\circ,k}(s),\bar Y^{\circ,k}(s),\bar R^{\circ,k}(s),\bar \Psi^{\circ,k}_1(s),\bar \Psi^{\circ,k}_2(s),\hat\Phi^k_1(s),\hat\Phi^k_2(s):0\le s\le t\}$, and
\begin{align}\label{new41}
\bar\Psi^{\circ,k}_1=\sigma^{-1}\sum_{i=1}^I  \theta_i\la_i^{1/2}\hat\Phi^k_{1,i},\qquad
\bar\Psi^{\circ,k}_2=\sigma^{-1}\sum_{i=1}^I  \theta_i\mu_i^{1/2}\hat\Phi^k_{2,i}.
\end{align}
Moreover, there are processes $\{\hat\phi^k_{j,i}\}_{j,i}$, progressively measurable w.r.t.~$\calG^n_t$, such that
\begin{align}\label{new30}
\hat\Phi^k_{1,i}(\cdot)=\int_0^\cdot\hat\phi^k_{1,i}(t)dt, \qquad \hat\Phi^k_{2,i}(\cdot)=\int_0^\cdot\hat\phi^k_{2,i}(t)\rho_idt,
\end{align}
which satisfy $\max_{j,i}\|\hat\phi^k_{j,i}\|_\iy\le k$. 
%
Finally, $\hat\Q$-a.s.,
\begin{align}\label{new31}
\begin{split}\liminf_{n\to\iy}\int_0^\iy e^{-\varrho t}(\hat \psi^{n,k}_{1,i}(t))^2dt
&\ge \int_0^\iy e^{-\varrho t}(\hat \phi^{k}_{1,i}(t))^2dt,\\
\liminf_{n\to\iy}\int_0^\iy e^{-\varrho t}(\hat \psi^{n,k}_{2,i}(t))^2dT^n_i(t)
&\ge \int_0^\iy e^{-\varrho t}(\hat \phi^{k}_{2,i}(t))^2\rho_idt.
\end{split}
\end{align}
\end{lemma}
{\bf Proof.}
The Skorokhod mapping formulation of the pre-limit processes, the limit in \eqref{new35}, and the $\calC$-tightness are private cases of Lemma \ref{lem46}. This results imply the Skorokhod formulation of the limiting process in \eqref{new35a}. 

The expressions in \eqref{new41} follow since $\theta^n\to\theta$ and from \eqref{new42},
\begin{align}
\label{new44}
\Psi^{\circ,n,k}_1(\cdot):=\sigma^{-1}\sum_{i=1}^I\theta^n_i\la_i^{1/2}\hat\Psi^{n,k}_{1,i}(\tau^{n,k}\wedge\cdot)
,\qquad\Psi^{\circ,n,k}_2(\cdot):=\sigma^{-1}\sum_{i=1}^I\theta^n_i\mu_i^{1/2}\hat\Psi^{n,k}_{2,i}(\tau^{n,k}\wedge\cdot).
\end{align}%
By definition, for every $i\in[I]$ and $j\in\{1,2\}$, $|\hat\psi^{n,k}_{j,i}|\le k$ and therefore,
$\hat\Q$-a.s., $ \{\hat\Psi^{n,k}_{j,i}\}_{j,i,n}$ are all Lipschitz-continuous with the same Lipschitz constant $k$. Since Lipschitz-continuity implies absolutely continuity, from Section IV.17 of \cite{del-mey1978} we get the existence of progressively measurable processes $\{\hat\phi^k_{j,i}\}_{j,i}$ such that \eqref{new30} holds. 

The bounds in \eqref{new31} follow immediately from Lemma \cite[Lemma A.3]{ASAFCOHEN2018}, noticing that by Skorokhod representation theorem the convergence $(T^{n,k}_i,\hat\Psi^{n,k}_{2,i})\To(\bar T_i,\hat\Phi^k_{2,i})$ can be replaced by u.o.c.~$\hat\Q$-a.s.~convergence.

Finally, we prove that $\bar W^{\circ,k}$ is an $(m,\sigma)$-Brownian motion. The martingale central limit theorem (\cite[Theorem 7.1.4]{ethkur}) implies that $\{(\check A^{n,k},\check  S^{n,k})\}_n$ converges to a $2I$-dimensional $\left(0,\tilde{\sigma}\right)$-Brownian motion, where $\check S^n=(\check S^n_i:i\in[I])$ and 
\begin{align}
\notag
\check S^{n,k}_i(t) :=n^{-1/2}\left(S^n_i(t)-\int_0^t\psi^n_{2,i}(s)ds\right),\quad\tilde{\sigma}=\text{Diag}(\la_1^{1/2},\ldots,\la_I^{1/2},\mu_1^{1/2},\ldots,\mu_I^{1/2}). 
\end{align}
Using a lemma regarding random change of time from \cite[p.151]{Bill} with $T^n_i\to\bar T_i$, gives that 
$\{(\check A^{n,k},\check D^{n,k})\}$ converges to a $(0,\check\sigma)$-Brownian motion, where 
\begin{align}
\notag
\check\sigma:=(\la^{1/2}_1,\ldots,\la^{1/2}_I,\la^{1/2}_1,\ldots,\la^{1/2}_I).
\end{align}
As a result from the definitions of $\hat W^n$, $W^{\sharp,n,k}$, and $W^{\circ,n,k}$ (see their equivalences in \eqref{newnew16}, \eqref{newnew16b}, and \eqref{newnew16c}), we get that $\bar W^{\circ,k}$ is an $(m,\sigma)$-Brownian motion w.r.t.~its own filtration. The proof that its filtration can be replaced by $\calF_t$ follows by the same lines of Proposition \ref{prop43} and therefore omitted.

%


\qed

The next proposition tells us that by truncating the intensities, the expected cost does not change much.
\begin{proposition}\label{prop41}
The following limit holds
\begin{align}\label{newnew9}
\begin{split}
&\limk\limsup_{n\to\iy}\left| \E^{\hat\Q}\Big[
\int_0^\iy e^{-\varrho t}\{h^a(X^{\circ,n}(t))+\varrho r R^{\circ,n}(t)
 \}dt
 \right.\\
 &\qquad\qquad\qquad
 \left.-
\int_0^\iy e^{-\varrho t}\{h^a(X^{\circ,n,k}(t))+\varrho r R^{\circ,n,k}(t)
 \}dt\Big]
\right|=0
\end{split}
\end{align}
and for every $i\in[I]$,
\begin{align}\label{newnew17}
\limk\limsup_{n\to\iy} \left(\left|L^\varrho_1(\hat \Q^{n}_{1,i}\|\PP^n_{1,i})-L^\varrho_1(\hat \Q^{n,k}_{1,i}\|\PP^n_{1,i})\right|+\left|L^\varrho_2(\hat \Q^{n}_{2,i}\|\PP^n_{1,i})-L^\varrho_2(\hat \Q^{n,k}_{2,i}\|\PP^n_{1,i})\right|\right)=0,
\end{align}
where $\hat\Q^{n,k}_{1,i}$ and $\hat\Q^{n,k}_{2,i}$ are the measures associated with $\psi^n_{1,i}(t)=\la^n_it+(\la_i n)^{1/2}\hat\psi^{n,k}_{1,i}(t)$ and $\psi^n_{2,i}(t)=\mu^n_it+(\mu_i n)^{1/2}\hat\psi^{n,k}_{2,i}(t)$. 
\end{proposition}
{\bf Proof.}
From the representations in \eqref{new32} and \eqref{new33}, the limits \eqref{new34} and \eqref{new35}, the definitions \eqref{new42}, \eqref{newnew16c}, and \eqref{new44}, and Lemma \ref{lem_Skorokhod}  it follows that in order to obtain \eqref{newnew9}, it is sufficient to show that
\begin{align}\label{newnew5}
&\limk\limsup_{n\to\iy}\E^{\hat\Q}\Big[\int_0^\iy e^{-\varrho t}\|W^{\circ,n,k}-W^{\circ,n}\|_tdt\Big]=0,\\\label{newnew18}
&\limk\limsup_{n\to\iy}\E^{\hat\Q}\Big[\int_0^\iy e^{-\varrho t}\|\hat\Psi^{n,k}_{1,i}-\hat\Psi^n_{1,i}\|_tdt\Big]=0,\\\label{newnew19}
&\limk\limsup_{n\to\iy}\E^{\hat\Q}\Big[\int_0^\iy e^{-\varrho t}\|\hat\Psi^{n,k}_{2,i}-\hat\Psi^n_{2,i}\|_tdT^n_i(t)\Big]=0,
\end{align}
Set the functions
$f^n_{1,i}:\left(-\la^n_i(\la_in)^{-1/2},\iy\right)\to\R$ and $f^n_{2,i}:\left(-\mu^n_i(\mu_in)^{-1/2},\iy\right)\to\R$ given by
\begin{align}
\notag
f^n_{1,i}(x):=\left(\la^n_i+(\la_in)^{1/2}x\right)\log\left(1+(\la_in)^{1/2}x/\la^n_i)\right)
-(\la_in)^{1/2}x,\\\notag
f^n_{2,i}(x):=\left(\mu^n_i+(\mu_in)^{1/2}x\right)\log\left(1+(\mu_in)^{1/2}x/\mu^n_i)\right)
-(\mu_in)^{1/2}x.
\end{align}
To obtain \eqref{newnew17}, recall \eqref{newnew2}. Then it is sufficient to show that 
\begin{align}\label{newnew20b}
&\limk\limsup_{n\to\iy}\E^{\hat\Q}\Big[\int_0^\iy e^{-\varrho t}|f^n_{1,i}(\hat\psi^{n,k}_{1,i}(t))-f^n_{1,i}(\hat\psi^n_{1,i}(t))|dt\Big]=0,\\\label{newnew20a}
&\limk\limsup_{n\to\iy}\E^{\hat\Q}\Big[\int_0^\iy e^{-\varrho t}|f^n_{2,i}(\hat\psi^{n,k}_{2,i}(t))-f^n_{2,i}(\hat\psi^n_{2,i}(t))|dT^n_i(t)\Big]=0,
\end{align}
We start with the limit in \eqref{newnew5}
Simple algebraic manipulation gives the bound
\begin{align}\notag
\|W^{\circ,n,k}-W^{\circ,n}\|_t\le C(\|\check A^{n,k}-\check A^n\|_t+\|\check D^{n,k}-\check D^n\|_t),
\end{align}
where $C$ refers to a finite positive constant
that is independent of $n$ and $t$ and which can change from one line to the next. From \eqref{new36} and Burkholder--Davis--Gundy inequality we get that for every $i\in[I]$,
\begin{align}\notag
\E^{\hat\Q}[\|\check A^{n,k}_i-\check A^n_i\|^2_t]
\le
Cn^{-1/2}
\E^{\hat\Q}\Big[\int_0^t|\hat\psi^{n}_{1,i}(s)|\one_{\{|\hat\psi^n_{1,i}(s)|>k\}}ds\Big].
\end{align}
Now, changing the order of integration gives
\begin{align}\label{newnew26}
&\limk\limsup_{n\to\iy}\E^{\hat\Q}\Big[\int_0^\iy e^{-\varrho t}\|\check A^{n,k}_i-\check A^n_i\|^2_tdt\Big]\\\notag
&\qquad\le
\limk\limsup_{n\to\iy}Cn^{-1/2}\E^{\hat\Q}\Big[\int_0^\iy e^{-\varrho t}
\Big\{\int_0^t|\hat\psi^{n}_{1,i}(s)|\one_{\{|\hat\psi^n_{1,i}(s)|>k\}}ds\Big\}
dt\Big]\\\notag
&\qquad
\le
\limk\limsup_{n\to\iy}Cn^{-1/2}\E^{\hat\Q}\Big[\int_0^\iy \varrho e^{-\varrho t}
|\hat\psi^{n}_{1,i}(t)|\one_{\{|\hat\psi^n_{1,i}(t)|>k\}}
dt\Big]=0.
\end{align}
Pay attention that by Proposition \ref{lem44}, the last limit holds even without the $n^{-1/2}$ term.

 The rest of the sufficient limits, namely \eqref{newnew18}, \eqref{newnew19}, \eqref{newnew20b}, and \eqref{newnew20a}, are treated similarly: notice that
 \begin{align}\notag
\|\hat\Psi^{n,k}_{1,i}-\hat\Psi^n_{1,i}\|_t
&\le C\sum_{i=1}^I\int_0^t|\hat\psi^{n}_{1,i}(s)|\one_{\{|\hat\psi^n_{1,i}(s)|>k\}}ds,\\\notag
|f^n_{1,i}(\hat\psi^{n,k}_{1,i}(t))-f^n_{1,i}(\hat\psi^n_{1,i}(t))|
&\le C\sum_{i=1}^If^n_{1,i}(\hat\psi^n_{1,i}(t))\one_{\{|f^n_{1,i}(\hat\psi^n_{1,i}(t))|>k\}},
\end{align}
and similar bounds apply for $j=2$ as well. We may continue in the same way as we did in \eqref{newnew26}, where now the $n^{1-2}$ term is absent, using the two bounds in Proposition \ref{lem44}.

\qed

We now show that the change of measure penalty in the truncated case can be approximated by a quadratic form. From \eqref{newnew2}, the bound $\sup_{t\in\R_+}|\psi^{n,k}_{1,i}(t)|\le k$, and Taylor's expansion of $\log(1+x)$, one has 
\begin{align}\notag
L^\varrho_1(\hat\Q^{n,k}_{1,i}\|\hat\PP^n_{1,i})
&=\E^{\hat\Q}\Big[\frac{1}{2}\int_0^\iy e^{-\varrho t}(\hat\psi^{n,k}_{1,i}(t))^2dt\Big]+o(1)
.
\end{align}
Similarly,
\begin{align}\notag
L^\varrho_2(\hat\Q^{n,k}_{2,i}\|\hat\PP^n_{2,i})&=
\E^{\hat\Q}\Big[\frac{1}{2}\int_0^\iy e^{-\varrho t}(\hat\psi^{n,k}_{2,i}(t))^2dT^n_i(t)\Big]+o(1)
.
\end{align}
From \eqref{new31} it follows that,
 \begin{align}\label{newnew7}
&\liminf_{n\to\iy}\; \Big[\sum_{i=1}^I\frac{1}{\kappa_{1,i}}L^\varrho_1(\hat\Q^{n,k}_{1,i}\|\hat\PP^n_{1,i})+\sum_{i=1}^I\frac{1}{\kappa_{2,i}}L^\varrho_2(\hat\Q^{n,k}_{2,i}\|\hat\PP^n_{2,i})\Big]\\\notag
&\qquad\ge
\sum_{i=1}^I\frac{1}{\kappa_{1,i}}
\E^{\hat\Q}\Big[\frac{1}{2}\int_0^\iy e^{-\varrho t}(\hat\phi^k_{1,i}(t))^2dt\Big]
+
\sum_{i=1}^I\frac{1}{\kappa_{2,i}}
\E^{\hat\Q}\Big[\frac{1}{2}\int_0^\iy e^{-\varrho t}(\hat\phi^k_{2,i}(t))^2\rho_idt\Big]+o(1)
.\end{align}
The identity
\begin{align}\notag
\frac{1}{k_1}\al_1^2+\frac{1}{k_2}\al_2^2\ge \frac{1}{k_1+k_2}(\al_1-\al_2)^2, \qquad k_1,k_2>0,\quad\al_1,\al_2\in\R,
\end{align}
implies that
\begin{align}\label{newnew8}
\frac{1}{2\kappa_{1,i}} (\hat\phi_{1,i}^k(t))^2+\frac{1}{2\kappa_{2,i}}\rho_i (\hat\phi_{2,i}^k(t))^2
\ge
\frac{1}{4\kaboom_i}\left(\hat\phi_{1,i}^k(t)-\rho_i^{1/2}\hat\phi_{2,i}^k(t)\right)^2,
\end{align}
where recall that  $\kaboom_i=\frac{1}{2}(\kappa_{1,i}+\kappa_{2,i})$.
Set
$$ \phi^{\sharp,k}:=\sigma^{-1}\sum_{i=1}^I  \theta_i\left(\la_i^{1/2}\hat\phi^k_{1,i}- \mu_i^{1/2}\hat\phi^k_{2,i}\rho_i\right).$$
By Cauchy--Schwartz inequality,
\begin{align}\notag
&\Big[\sigma^{-2}\sum_{i=1}^I(\theta\hat\sigma)_i^2\kaboom_i\Big]\times\sum_{i=1}^I\frac{1}{2\kaboom_i}\left(\hat\phi_{1,i}^k(t)-\rho_i^{1/2}\hat\phi_{2,i}^k(t)\right)^2
\ge
(\phi^{\sharp,k}(t))^2.
\end{align}
Then,
\begin{align}\notag
\sum_{i=1}^I\frac{1}{4\kaboom_i}\left(\hat\phi_{1,i}^k(t)-\rho_i^{1/2}\hat\phi_{2,i}^k(t)\right)^2
\ge
\frac{1}{2\eps}(\phi^{\sharp,k}(t))^2,
\end{align}
where recall that  $\eps=\sigma^{-2}\sum_{i=1}^I(\theta\hat\sigma)_i^2\kaboom_i$.

From \eqref{newnew7}, \eqref{newnew8}, and the last bound, we obtain that
\begin{align}\label{newnew20}
\begin{split}
&\liminf_{n\to\iy}\;\Big[ \sum_{i=1}^I\frac{1}{\kappa_{1,i}}L^\varrho_1(\hat\Q^{n,k}_{1,i}\|\hat\PP^n_{1,i})+\sum_{i=1}^I\frac{1}{\kappa_{2,i}}L^\varrho_2(\hat\Q^{n,k}_{2,i}\|\hat\PP^n_{2,i})\Big]\\
&\qquad\ge \frac{1}{2\eps}\E^{\hat \Q}\Big[\int_0^\iy e^{-\varrho t}(\phi^{\sharp,k}(t))^2dt\Big].
\end{split}
\end{align}

Fix $\delta_1>0$. Combining the last limit with the ones from \eqref{new40} and \eqref{newnew9}
 and recalling \eqref{new26} and the bound $\|R^{\circ,n}\|\le\|\theta^n\|\|\hat R^n\|$, we get that there is $k_{\delta_1}>0$ such that for every $k\ge k_{\delta_1}$,
\begin{align}\label{newnew10}
 %
 &\limsup_{n\to\iy}\;J^n( X^n(0),U^n(a),R^n(a),\hat\Q^n;\kappa)\\\notag
 &\quad\le
 \E^{\hat\Q}\Big[
\int_0^\iy e^{-\varrho t}\{h^a(\bar X^{\circ,k}(t))+\varrho r\bar R^{\circ,k}(t)
 -\frac{1}{2\eps}(\phi^{\sharp,k}(s))^2\}dt\Big]+\delta_1/2,
 \end{align}
where notice that \eqref{new35a}--\eqref{new30}, and the definition of $\phi^{\sharp,k}$ gives
\begin{align}\notag
(\bar X^{\circ,k},\bar Y^{\circ,k},\bar R^{\circ,k})=\Gamma_{[0,\aey ]}\Big(x_0+\bar W^{\circ,k}+\sigma\int_0^\cdot\phi^{\sharp,k}(s)ds\Big).
\end{align}
Fix such $k$.
Since $\aey \to \Barr$ as $\delta_0\to0$ (see the paragraph before \eqref{eq2304}), Lemma \ref{lem_Skorokhod}
and \eqref{eq2307b} imply that
\begin{align}\notag
&\int_0^\iy e^{-\varrho t}\{h^a(\bar X^{\circ,k}(t))+\varrho r\bar R^{\circ,k}(t)
 -\frac{1}{2\eps}(\phi^{\sharp,k}(s))^2\}dt\\\notag
 &\qquad
 \overset{\hat\Q}{\Longrightarrow}
 \int_0^\iy e^{-\varrho t}\{h(\bar X(t))+\varrho r\bar R(t)
 -\frac{1}{2\eps}(\phi^{\sharp,k}(s))^2\}dt,
\end{align}
as  $\delta_0\to0$, where
\begin{align}\notag
(\bar X,\bar Y,\bar R)=\Gamma_{[0,\Barr]}\Big(x_0+\bar W^{\circ,k}+\sigma\int_0^\cdot\phi^{\sharp,k}(s)ds\Big).
\end{align}
{\bf Obtaining the upper bound.} 
As argued earlier, to conclude the convergence in expectation of the above it is sufficient to show that
\begin{align}\label{newnew4}
\limsup_{\aey \to \Barr}\E\Big[\int_0^\iy e^{-\varrho t}\|\bar R^{\circ,k}(t)\|^2dt\Big]<\iy.
\end{align}
The proof is borrowed from \cite[(85)]{ata-shi} and adapted to our case with the additional stochastic drift $\sigma\int_0^\cdot\phi^{\sharp,k}(s)ds$.
Apply It\^{o}'s formula to $(\bar X^{\circ,k})^2$, use the reflection conditions $\int_0^t\bar X^{\circ,k}(s)d\bar Y^{\circ,k}(s)=0$
and $\int_0^t\bar X^{\circ,k}d\bar R^{\circ,k}(s)=\aey \bar R^{\circ,k}(t)$, to get
\[
\bar R^{\circ,k}(t)=\frac{1}{2\aey }\Big\{(\bar X^{\circ,k}(0 ))^2-(\bar X^{\circ,k}(t))^2+2\int_0^t\bar X^{\circ,k}(s)[d\bar W^{\circ,k}(s)+\sigma\phi^{\sharp,k}(s)ds]
+\sigma^2t\Big\}.
\]
Since $\bar X^{\circ,k}$ and $
\sigma\phi^{\sharp,k}(s)$ are bounded, 
\eqref{newnew4} follows easily.

Recall the definition of $J$ and $V$ from Section \ref{sec32} and that $\bar X^{\circ,k}(0)=x_0$, then
\begin{align}\notag
&\lim_{\delta_0\to 0}\E^{\hat\Q}\Big[\int_0^\iy e^{-\varrho t}\{h^a(\bar X^{\circ,k}(t))+\varrho r\bar R^{\circ,k}(t)
 -\frac{1}{2\eps}(\phi^{\sharp,k}(t))^2\}dt\Big]\\\notag
 &\qquad
=
\E^{\hat\Q}\Big[\int_0^\iy e^{-\varrho t}\{h(\bar X(t))+\varrho r\bar R(t)
 -\frac{1}{2\eps}(\phi^{\sharp,k}(t))^2\}dt\Big]\\\notag
 &\qquad
= J(\bar X^{\circ,k}(0), \bar Y,\bar R,\Q;\eps)\\\notag
&\qquad \le V(x_0;\eps),
\end{align}
where $\Q$ is the measure associated with the rate $\phi^{\sharp,k}$. The last inequality follows by the optimality of the minimizer's reflected strategy, see Proposition \ref{prop31}. 

Together with \eqref{newnew10},
we obtain that for sufficiently small $\delta_0>0$ and large $n$,
\begin{align}\notag
 %
 &J^n( X^n(0),U^n(a),R^n(a),\hat\Q^n;\kappa)\le V( x_0;\eps)+\delta_1.
 \end{align}
 Since the constant $\delta_1(>0)$ and the measures $\{\hat\Q^n\}_n$ were chosen arbitrary
 , we get \eqref{new5}.


\subsection{The maximizer's asymptotic behavior}\label{sec44}
Recall that under the measure $\hat\Q^n_{1,i}$ (resp., $\hat\Q^n_{2,i}$), $A^n_i$ (resp., $S^n_i(T^n_i)$) is a counting process with intensity $\psi^n_{1,i}$ (resp., $\psi^n_{2,i}U^n_i$). Notice also that under the measures $\hat\Q^n_{j,i}$, $j\in\{1,2\}$, the critically load condition might be violated since we do not restrict the intensities $\{\hat\psi^n_{j,i}\}_{j,i,n}$ in such a way. However, as follows from Sections \ref{sec42} and \ref{sec43}, such changes of measures are `too costly' and will be avoided by the maximizer so that, {\bf as a consequence}, `in average', 
$\psi^n_{1,i}(t)=\la_i^n+\calO(n^{1/2})$ and $
\psi^n_{2,i}(t)=\mu_i^n+\calO(n^{1/2}).$ 
Indeed, in the proof of the lower bound, we used the equilibrium control for the maximizer from the RSDG and set up the maximizer's control in the QCP to be given by \eqref{501}--\eqref{newnew3}. Therefore, the critically load condition holds. 

In the proof of the upper bound, we show in Section \ref{sec431} and more specifically in Proposition \ref{lem44} that 
$\psi^n_{1,i}(t)=\la_i^n+(\la_in)^{1/2}\hat\psi^n_{i,1}(t)$, where the expectations $\E^{\hat\Q}\Big[\int_0^\iy e^{-\varrho t}
 |\hat\psi^n_{1,i}(t)|dt\Big]$ are uniformly bounded in $n$, and similarly for $\psi^n_{2,i}$.

The main reasoning for this phenomenon lies in the fact that we are studying a diffusion scaled problem and as such the uncertainty is expressed on the scaled system. Now, by observing the diffusion scaled system the DM can rule out measures under-which $|\psi^n_{i,1}-\la^n_i|+|\psi^n_{i,2}-\mu^n_i|$ is of higher order than $n^{1/2}$ (for sufficiently large $n$). This is because the diffusion coefficient in the limiting dynamics under the new measure would be different from the one that emerges from the reference measure and in the limiting problem this could be identified `very fast'. Another reason is of course the structure of the penalty using the Kullback--Leibler divergence, which also leads to a limiting problem with linear drift and quadratic penalty in the maximizer's control
, see \eqref{321}. 

\section{Future outlook}\label{sec5}
\beginsec
\begin{enumerate}
\item Recall Remark \ref{rem_21}.(ii). Much as in the Markovian case,  the heavy traffic approximation for G/G/1 also  scales \eqref{204a} to \eqref{new18}, whose form uses the {\rm mean} arrival and service rates $\la^n_i$ and $\mu^n_i$. In the change of measure given by 
\eqref{216} appear the {\rm actual} arrival and service rates under the measures $\PP^n_{j,i}$ and $\hat\Q^n_{j,i}$, namely $(\la^n_i,\mu^n_i)$ and $(\psi^n_{1,i}(\cdot),\psi^n_{2,i}(\cdot))$, respectively. 
In the Markovian case (M/M/1) the {\rm mean} and the {\rm actual} arrival and service rates are identical. Our technique heavily uses this identification. It is desirable to extend the main result of this paper beyond the Markovian setting, to
general service time distributions and renewal arrival distributions, under second moment
conditions. 

\item Our model considers finite buffers and linear cost. The author is currently testing the robustness of the {\it generalized} $c\mu$ policy in the unconstrained case. I expect to benefit from the pathwise minimality property of the Skorokhod mapping (see \cite{ata-bis}) in three aspects: (1) replacing the Kullback--Leibler divergence with a more general divergence satisfying basic properties, such as: growth, continuity, and convexity; (2) attaining the lower bound without needing the smoothness of the value function. This allows us to skip the differential equation analysis done in \cite{Cohen2017}; and (3) proving the lower bound without the time rescaling procedure.

The state-space collapse part in the proof of the upper bound in this case should be simpler again since rejections are not allowed. However, using a general divergence leads to a more subtle analysis of the maximizer's penalty part of the cost.

\item The present model assumes the uncertainty in the diffusion scale, and as discussed in Section \ref{sec44} the system stays critically loaded. The author studies a similar type of uncertainty in the fluid scale. In this case we do not assume a reference model, but rather the maximizer chooses the worst set of parameters (possibly stochastic and time-dependent) within a given set of possible parameters, without a penalty component. In this case it is less trivial that an index policy is asymptotically optimal.

\item Recently, \cite{blanchet2016quantifying} quantified the impact of model misspecification when computing general expected values using an optimal transport cost instead of divergences. It would be interesting to tackle the QCP presented here and challenging the $c\mu$ rule using optimal transport tools.

\end{enumerate}

\section{Proofs}\label{sec6}
\beginsec

{\bf Proof of Lemma \ref{lem41}.}
The $\cal C$-tightness of $\{\tilde \tau^n,\tilde Y^n,\tilde R^n)\}_n$ follows by the following observation, which relies on the estimate $|\tau^n(\tilde\tau^n(t))-t|\le\|\theta^n\|/\sqrt{n}$ and the definition of $\tau^n$.
\begin{align}\label{515}
\frac{2\|\theta^n\|}{\sqrt{n}}+t-s&\ge \tau^n(\tilde\tau^n(t))-\tau^n(\tilde\tau^n(s))\\\notag
&=\tilde\tau^n(t)-\tilde\tau^n(s)+\theta^n\cdot\tilde R^n(t)-\theta^n\cdot\tilde R^n(s)+\theta^n\cdot\tilde Y^n(t)-\theta^n\cdot\tilde Y^n(s).
\end{align}

Next, since the sequence $\{(\tilde\tau^n,\tilde Y^n)\}_n$ is $\cal C$-tight, we get by \eqref{207} and
the limit $\mu^n_in^{-1/2}\to\iy$, that $ \tilde T^n\overset{\hat \Q^n}\Longrightarrow \tilde T$.

From \eqref{505}, the $\calC$-tightness of $\{\tilde \tau^n,\tilde Y^n,\tilde R^n)\}_n$, and \eqref{new8}, the $\cal C$-tightness of $\{\tilde X^n\}_n$ follows once we show that $\{(\tilde A^n,\tilde D^n)\}_n$ is $\cal C$-tight. Recalling that $\{\tilde\tau^n\}_n$ is $\cal C$-tight, then in order to prove the latter statement, it is sufficient to show the $\cal C$-tightness of $\{(\check A^n,\check D^n)\}$. 
%
%

For every $i\in[I]$ and $n\in\N$, the processes $\psi^n_{1,i}$ and $\psi^n_{2,i}U^n_i$ are the intensities of $A^n_i$ and $D^n_i:=S^n_i(T^n_i)$, respectively.
Denote 
\begin{align}\notag
(W^n_1,\ldots,W^n_{2I})&:=(A^n_1,\ldots,A^n_I,D^n_1,\ldots,D^n_I),\\\notag
\check W^n=(\check W^n_1,\ldots,\check W^n_{2I})&:=(\check A^n_1,\ldots,\check A^n_I,\check D^n_1,\ldots,\check D^n_I),\\\notag
\tilde W^n=(\tilde W^n_1,\ldots,\tilde W^n_{2I})&:=(\tilde A^n_1,\ldots,\tilde A^n_I,\tilde D^n_1,\ldots,\tilde D^n_I).
\end{align}
Therefore, 
for every $j\in\{1,\ldots,2I\}$, $\check W^n_j$ is a martingale w.r.t.~its own filtration, under the measure $\hat\Q^n$. 
Notice that the quadratic variation of $\check W^n_j$ satisfies,
\begin{align}\label{517}
[\check W^n_j]=\frac{1}{n}W^n_j,
\end{align}
which is of order $t$ thanks to \eqref{newnew3}, \eqref{new8}, and \eqref{205}.
%
Fix an arbitrary $T>0$ and a stopping time $\pi$. Then, by
Burkholder--Davis--Gundy inequality (see \cite[Theorem 48]{Protter2004}), and \eqref{517}, we get,
\begin{align}\notag
\left(\E^{\hat\Q^n}\|\check W^n_j(\pi+\cdot)-\check W^n_j(\pi)\|_\delta\right)^2
&\le
\E^{\hat\Q^n}\|\check W^n_j(\pi+\cdot)-\check W^n_j(\pi)\|_\delta^2\le
C_3\E^{\hat\Q^n}[\check W^n_j(\pi+\cdot)]_\delta\le C_4\delta
\end{align}
for some constants $C_3,C_4>0$, independent of $n$, $\delta$, and $\pi$.
Therefore, Aldous criterion for tightness holds (see e.g., \cite[Theorem 16.10]{Bill}). Since the jumps of these processes are of order $O(n^{-1/2})$, any limit process
has continuous paths with probability $1$ and $\calC$-tightness of $\{\check W^n\}_n$ is proved. 

The first two identities in \eqref{514b} follow by the convergence $(\check W^n,\tilde W^n,\tilde\tau^n)
\Rightarrow (\check W,\tilde W,\tilde\tau)$, which in turn follows by the tightness of $\{(\check W^n,\tilde W^n,\tilde\tau^n)\}_n$. Finally, the quadratic variation of $\tilde B$ follows by \eqref{newnew3}, \eqref{514b}, and the martingale central limit theorem, see \cite[Theorem 7.1.4]{ethkur}. 

As mentioned after Lemma \ref{lem41}, by Skorokhod's representation theorem, we may assume that along a subsequence \eqref{521} holds a.s., u.o.c.~under some probability space that supports all the relevant processes. 
Now, from  Lemma \ref{lem_A2}, \eqref{501}--\eqref{newnew3}, and \eqref{505}, we obtain \eqref{514}.

\hfill$\Box$

{\bf Proof of Lemma \ref{lem42}.}
Fix $0<s<t$. From \eqref{511},
it follows that
\begin{align}\label{523}
\begin{split}
\Q^*(\tilde\tau^n(t)\le s)&=\Q^*(\tau^n(s)\ge t)=\Q^*(s+\theta^n\cdot\hat R^n(s)+\theta^n\cdot\hat Y^n(s)>t)\\
&\le \frac{1}{t-s}\E^{\Q^*}[\theta^n\cdot\hat R^n(s)+\theta^n\cdot\hat Y^n(s)].
\end{split}
\end{align}
From \eqref{new8}, \eqref{505}, and the inequality $T^n(u)\le u$, $u\in\R_+$ we get that there is a constant $C_5>0$ independent of $n$ and $t$, and such that 
\begin{align}\notag
\hat Y^n_i(s)\le C_5(1+\hat R^n_i(s)+\check A^n_i(s)-\check D^n_i(s)).
\end{align}
Recalling that $\check A^n_i-\check D^n_i$ is a martingale, and \eqref{504}, we get that the last display in \eqref{523} is bounded above by $\frac{C_6}{t-s}$, for some $C_6>0$, independent of $n$ and $t$.
Now, since for every $s\in\R_+$ and $0\le t\le t'$, $\{\tilde\tau(t)>s\}\supseteq \{\tilde\tau(t')>s\}$, we get that
\begin{align}\notag
\Q^*\left(\lim_{t\to\iy}\tilde\tau(t)<s\right)=\lim_{t\to\iy}\Q^*(\tilde\tau(t)<s).
\end{align}
Using the convergence in law $\tilde\tau^n\overset{d}{\Rightarrow}\tilde\tau$, and since $\Q^*(\tilde\tau^n(t)<s)\le\frac{C_6}{t-s}$, we conclude that
\begin{align}\notag
\Q^*\left(\lim_{t\to\iy}\tilde\tau(t)<s\right)\le\lim_{t\to\iy}\limsup_{n\to\iy}\Q^*(\tilde\tau^n(t)<s)=0.
\end{align}

\hfill$\Box$

{\bf Proof of Proposition \ref{prop43}.} 
We use the L\'{e}vy characterization for Brownian motions. Namely, we show that
$B^*(t)$ and $B^*(t)(B^*(t))^\top- \calI t$ are continuous local martingales w.r.t.~$\calG^*$, where recall that $\calI$ is the identity matrix of order $I\times I$.

To simplify the presentation, set
\begin{align}\notag
\check B^n&=(\check B^n_1,\ldots,\check B^n_I)=\hat\sigma^{-1}(\check A^n-\check D^n),\qquad\check B=(\check B_1,\ldots,\check B_I)=\hat\sigma^{-1}(\check A-\check D),\\\notag
\tilde B^n&=(\tilde B^n_1,\ldots,\tilde B^n_I)=\hat\sigma^{-1}(\tilde A^n-\tilde D^n),\qquad\tilde B=(\tilde B_1,\ldots,\tilde B_I)=\hat\sigma^{-1}(\tilde A-\tilde D).
\end{align}
We start by arguing the continuity.
From \eqref{521}, $(\check B^n,\tilde \tau^n,\tilde B^n)\to(\check B,\tilde\tau,\tilde B)$ $\Q^*$-a.s., u.o.c.~and therefore, $\tilde B(\cdot)=\check B(\tilde\tau(\cdot))$, $\Q^*$-a.s., u.o.c. Thus, $B^*(\cdot)=\tilde B(\tau(\cdot))=\check B(\tilde\tau(\tau(\cdot)))=\check B(\cdot)$, which is continuous $\Q^*$-a.s., see Lemma \ref{lem41}. The proof that $B^*(t)(B^*(t))^\top- \calI t$ is a local martingale follows by the same lines of the proof that  $B^*(t)$ is a local martingale. We start with the proof of the latter and then add the missing details for the former. Recall the definition of $\calG^n$ from Definition \ref{def21}.(i).

Fix $t,s\in\R_+$ and $n\in\N$. Notice that $\{\tilde\tau^n(s)\le t\}=\{\tau^n(t)\ge s\}=\{t+\theta\cdot\hat R^n(t)+\theta\cdot\hat Y^n(t)\ge s\in\calG^n_t$. Thus, $\tilde\tau^n(s)$ is a $\calG^n_t$-stopping time. Recall that $\check A^n$ and $\check D^n$ are $\calG^n_t$-martingales, then the optional sampling theorem (see Problem 1.3.24 in [20]) yields that
$\tilde B^n(t+\cdot)=\check B^n(\tilde\tau(t+\cdot))$ is a $\calG^n(\tilde\tau^n(t))$-martingale. As a consequence, for every $i\in[I]$ and every $\calG^n(\tilde\tau^n(t))$-measurable random variable $\zeta^n$,
\begin{align}\notag
\E^{\Q^*}[(\tilde B^n_i(t+s)-\tilde B^n_i(t))\zeta^n]=0.
\end{align}
Recall that $t\mapsto \tilde \tau^n(t)$ is nondecreasing and therefore $\calG^n(\tilde\tau^n(t))$ is a filtration. Now,
for every bounded continuous function $g$,
\begin{align}\notag
&g(\tilde A^n(s_m),\tilde D^n(s_m),\tilde Y^n(s_m),\tilde R^n(s_m),\tilde\tau^n(s_m);\;0\le s_m\le t, m=1,\ldots,k)\\\notag
&\quad=g(\check A^n(\tilde\tau^n(s_m)),\check D^n(\tilde\tau^n(s_m)),\hat Y(\tilde\tau^n(s_m)),\hat R(\tilde\tau^n(s_m)),\tilde\tau^n(s_m);\;0\le s_m\le t, m=1,\ldots,k)
\end{align}
 is $\calG^n(\tilde\tau^n(t))$-measurable, where $k\in\N$, and $\{s_m\}$ is an arbitrary sequence. Combining the last two displays one gets,
 \begin{align}\notag
\E^{\Q^*}[&g(\tilde A^n(s_m),\tilde D^n(s_m),\tilde Y^n(s_m),\tilde R^n(s_m),\tilde\tau^n(s_m);\;0\le s_m\le t, m=1,\ldots,k)\\\label{537}
&\qquad\qquad\qquad\qquad\qquad\qquad\qquad\qquad\qquad\qquad\qquad\times(\tilde B^n_i(t+s)-\tilde B^n_i(t))]=0.
\end{align}
From \eqref{505},
\begin{align}\label{534}
\tilde X^n_i(t)
&=\tilde X^n_i(0)+ \hat m^n_i \tilde\tau^n(t)+\tilde A^n_i(t)-\tilde D^n_i(t)+\tilde Y^n_i(t)-\tilde R^n_i(t)\\\notag
&\qquad+\la_i^{1/2}\int_0^{\tilde\tau^n(t)}\hat\psi^n_{1,i}(s)ds -\mu_i^{1/2}\int_0^{\tilde\tau^n(t)}\hat\psi^n_{2,i}(s)dT^n_i(s).
\end{align}
From \eqref{new8}, \eqref{534}, and the bound $T^n(u)\le u$, $u\in\R_+$, we get that there are constants $C_7,C_8>0$, independent of $t$, such that
for sufficiently large $n$ one has
\begin{align}\notag
|\tilde B^n_i(t)|\le C_7(\theta^n\cdot\tilde R^n_i(t)+\theta^n\cdot\tilde Y^n_i(t)+\tilde \tau^n_i(t)+1)\le C_8(t+1),\quad t\in\R_+,
\end{align}
where the last inequality follows by the same arguments leading to \eqref{515}.
From \eqref{521}, \eqref{537}, and the bounded convergence theorem, we get that,
\begin{align}\notag
\E^{\Q^*}[&g(\tilde A(s_m),\tilde D(s_m),\tilde Y(s_m),\tilde R(s_m),\tilde\tau(s_m);\;0\le s_m\le t, m=1,\ldots,k)\\\label{540}
&\qquad\qquad\qquad\qquad\qquad\qquad\qquad\qquad\qquad\qquad\qquad\times(\tilde B_i(t+s)-\tilde B_i(t))]=0,
\end{align}
which implies that $\tilde B$ is a $\tilde\calG_t$-martingale.

Next, for every $i\in[I]$, we show that the composition $B^*_i(t)=\tilde B_i(\tau(t))$ is a $\calG^*_t$-local martingale. For this, we need to define the following stopping times. Fix $M>0$ and set,
\begin{align}\notag
\tilde\pi_{i,M}&:=\inf\{t\ge 0: \tilde B_i(t)> M\},\\\notag
\pi_{i,M}&:=\inf\{t\ge 0: \check B_i(t)> M\}.
\end{align}
One can verify that $\tau(\pi_{i,M})=\tilde\pi_{i,M}$.

We now show that $B^*_i(t\wedge \pi_{i,M})$, which by definition equals $\tilde B_i(\tau(t\wedge \pi_{i,M}))$, is a $\calG^*_t$-martingale. Since $\lim_{M\to\iy}\pi_{i,M}=\iy$, $\Q^*$-a.s., we conclude that $B^*$ is a $\calG^*_t$-local martingale. For this, set,
$\tilde B_{i,M}(\cdot):=\tilde B_i(\cdot\wedge\tilde\pi_{i,M})$.
We use the optional sampling theorem given in \cite[Theorem 2.2.13]{ethkur}, which (in adaptation to our notation) states that if for every $t\in\R_+$,
\begin{align}\label{538}
\E^{\Q^*}|\tilde B_{i,M}(\tau(t\wedge \pi_{i,M}))|<\iy
\end{align}
and
\begin{align}\label{539}
\lim_{T\to\iy}\E^{\Q^*}\left[|\tilde B_{i,M}(T)|\one_{\{\tau(t)>T\}}\right]=0,
\end{align}
then for every $0\le s\le t$,
\begin{align}\notag
\E^{\Q^*}\left[\tilde B_{i,M}(\tau(t\wedge \pi_{i,M}))\mid\calG^*(s\wedge \pi_{i,M}))\right]&=\E^{\Q^*}\left[\tilde B_{i,M}(\tau(t\wedge \pi_{i,M}))\mid\tilde\calG(\tau(s\wedge \pi_{i,M})))\right]\\\notag
&=\tilde B_{i,M}(\tau(s\wedge \pi_{i,M})),
\end{align}
where we used the identity $\calG^*(t)=\tilde\calG(\tau(t))$. 
Therefore, $\tilde B_{i,M}(\tau(t\wedge \pi_{i,M}))$ is a $\calG^*_t$-martingale. Notice that
\begin{align}\notag
B^*_i(t\wedge \pi_{i,M})&=\tilde B_i(\tau(t\wedge\pi_{i,M}))=\tilde B_i(\tau(t)\wedge\tau(\pi_{i,M}))=
\tilde B_i(\tau(t)\wedge\tau(\pi_{i,M})\wedge \tilde\pi_{i,M})\\\notag
&=\tilde B_i(\tau(t\wedge\pi_{i,M}))\wedge \tilde\pi_{i,M})=\tilde B_{i,M}(\tau(t\wedge \pi_{i,M}))
\end{align}
and therefore, $B^*_i(t\wedge \pi_{i,M})$ is a $\calG^*_t$-martingale.
Indeed, the first and last equalities follow by the definitions of $B^*$ and $\tilde B_{i,M}$, respectively. The second and the forth equalities follow by the monotonicity of the function $t\mapsto\tau(t)$. Finally, the third equality follows since $\tau(\pi_{i,M})=\tilde\pi_{i,M}$.

We now prove that Properties \eqref{538} and \eqref{539} hold. Property \eqref{538} follows by the definition of $\tilde B_{i,M}$. 
To prove Property \eqref{539} notice that
\begin{align}\notag
\E^{\Q^*}\left[|\tilde B_{i,M}(T)|\one_{\{\tau(t)>T\}}\right]\le \E^{\Q^*}\left[|\tilde B_{i,M}(T)|\one_{\{\tilde\tau(T)\le t\}}\right]\le M\Q^*(\tilde\tau(T)\le t).
\end{align}
Now, from Lemma \ref{lem42} the leff-hand side of the above approaches $0$ as $T\to\iy$ and \eqref{539} is proven.

We end the proof by providing the missing arguments for the proof that $B^*(t)(B^*(t))^\top- \calI t$ is a martingale. One may go over the proof and replace the $\tilde B_i(t)$'s, $i\in[I]$ with $\tilde N_{ij}(t)=\tilde B_{i}(t)\tilde B_{j}(t)-\delta_{ij}t$,
$i,j\in[I]$. The only difference between the proofs lies in proving that $\tilde N(t)=\tilde B(t)(\tilde B(t))^\top-\calI t$ is a $\tilde\calG_t$-martingale. Or equivalently, in showing that \eqref{540} holds with $\tilde N_{ij}$
, replacing $\tilde B_i$. To this end, recall that $\hat\Q^n=\prod_{i=1}^I(\hat\Q^n_{1,i}\times\hat\Q^n_{2,i})$. Thus, $\{A^n_1,\ldots,A^n_{I},D^n_1,\ldots,D^n_{I}\}$ are mutually independent under $\hat Q^n$ and therefore also under $\Q^*$. Moreover, using the continuity of $\tilde\tau^n$ and the notation $\Delta L(t)=L(t)-L(t-)$, we get that for any $i,j\in[I]$, $i\ne j$, the following equalities hold $\Q^*\text{-a.s.}$,
\begin{align}\notag
[\tilde B^n_i,\tilde B^n_j](t)=\sum_{0\le s\le t}\Delta\tilde B^n_i(s)\Delta\tilde B^n_j (s)=\frac{1}{n}\sum_{0\le s\le t}\Delta B^n_i(\tilde\tau^n(s))\Delta B^n_j (\tilde\tau^n(s))=0
\end{align}
and
\begin{align}\notag
[\tilde B^n_i,\tilde B^n_i](t)&=[\tilde B^n_i](t)=\sum_{0\le s\le t}(\Delta\tilde B^n_i(s))^2=\frac{1}{n\hat\sigma_i^2}\sum_{0\le s\le t}(\Delta A^n_i(\tilde\tau^n(s)))^2+(\Delta D^n_i(\tilde\tau^n(s)))^2 \\\notag
&=\frac{1}{2n\la_i}\left(A^n_i(\tilde\tau^n(t))+D^n_i(\tilde\tau^n(t))\right)\underset{n\to\iy}{\longrightarrow}\tilde\tau(t),
\end{align}
where the limit holds by the strong law of large numbers.
Since $\tilde B^n_i(t)\tilde B^n_j(t)-[\tilde B^n_i,\tilde B^n_j](t)$ is a $\calG^n(\tilde\tau(t))$-martingale, by taking the limit $n\to\iy$, one deduces that \eqref{540} holds with $\tilde N_{ij}$ replacing $\tilde B_i$.

\hfill$\Box$

{\bf Proof of Lemma \ref{lem45b}.}
From Proposition \ref{lem44} it follows immediately that for every $T>0$,
\begin{align}\notag
\lim_{K\to\iy}\limn \hat\Q\left(\|(\hat\Psi^n_1,\hat\Psi^n_2))\|_T\ge K\right)=0.
\end{align}
Fix $\delta,\eta,K>0$ and for every $n\in\N$ set  \begin{align}
\label{newnew15}
P^n:=\int_0^te^{-\varrho t}\left((\hat\psi^n_{1,i}(s))^2ds+(\hat\psi^n_{2,i}(s))^2dT^n_i(s)\right).
\end{align}
Clearly, for any $K>0$,
\begin{align}\label{new15}
\begin{split}
\hat\Q\left(\osc_T((\hat\Psi^n_1,\hat\Psi^n_2),\delta)\ge \eta\right)
&=
\hat\Q\left(\osc_T((\hat\Psi^n_1,\hat\Psi^n_2),\delta)\ge \eta, P^n> K\right)\\
&\quad+\hat\Q\left(\osc_T((\hat\Psi^n_1,\hat\Psi^n_2),\delta)\ge \eta, P^n\le K\right)
.
\end{split}
\end{align}
Proposition \ref{lem44} implies that
\begin{align}\label{new16}
\begin{split}
&\lim_{K\to\iy}\limsup_{\delta\to0+}\limsup_{n\to\iy}\hat\Q\left(\osc_T((\hat\Psi^n_1,\hat\Psi^n_2),\delta)\ge \eta, P^n> K\right)\\
&\quad\le
\lim_{K\to\iy}\limsup_{n\to\iy}\hat\Q\left( P^n> K\right)=0.
\end{split}
\end{align}

We now examine the second term on the r.h.s.~of \eqref{new15}. On
the event $\{P^n\le K\}$, Jensen's inequality implies that there exists a constant $C_T>0$, independent of $n$, such that for every $0\le s<t\le T$,
\begin{align}\notag
\|(\hat\Psi^n_1,\hat\Psi^n_2)(t)-(\hat\Psi^n_1,\hat\Psi^n_2)(s)\|^2\le C_T(t-s)P^n\le C_TK(t-s).
\end{align}
As a result,
\begin{align}\notag
\lim_{\delta\to0+}\limn\hat\Q\left(\osc_T((\hat\Psi^n_1,\hat\Psi^n_2),\delta)\ge \eta, P^n\le K\right)=0.
\end{align}
Combining it with \eqref{new16}, we get
\begin{align}\notag
\lim_{\delta\to0+}\lim_{n\to\iy}\hat\Q\left(\osc_T((\hat\Psi^n_1,\hat\Psi^n_2),\delta)\ge \eta\right)=0,
\end{align}
and the $\calC$-tightness is established.

\qed

{\bf Proof of Lemma \ref{lem45}.}
By \eqref{newnew24} we have,
\begin{align}\label{new20}
\|\hat Y^n\|_T\le \|\hat X^n\|_T+\|\check A^n\|_T+\|\check D^n\|_T+\|\hat m^n\|T+C_2\|(\hat\Psi^n_1,\hat\Psi^n_2)\|_T+\|\hat R^n\|_T,
\end{align}
where $C_2$ depends solely on $\{(\la_i,\mu_i)\}_i$. 
From \eqref{new19} it follows that there exists a constant $C_3>0$, independent of $n$ and $T$, such that for every $T\in\R_+$,
\begin{align}\label{new21}
\|\hat R^n\|_T\le C_3\Big(1+T+\|\check A^n\|_T+\|\check  D^n\|_T +\|(\hat\Psi^n_1,\hat\Psi^n_2)\|_T \Big).
\end{align}
Recall that $\|\hat X^n\|_T$ is bounded and from Lemma \ref{lem45b}, $\|(\hat\Psi^n_1,\hat\Psi^n_2)\|_T$ is tight. Then, once we show that $\{(\check A^n,\check D^n)\}_n$ is $\calC$-tight, from \eqref{new20} and \eqref{new21}, we get that
$\{\|\hat Y^n\|_T\}_n$ is tight. 
Now, by the definitions of $\hat Y^n$ and $\mu^n$, see \eqref{207} and \eqref{205}, we get that $\{T^n_i\}_n$ converges u.o.c. to $\bar T_i$, where $\bar T_i(t)=\rho_it$.

We now argue the $\calC$-tightness of $\{(\check A^n,\check D^n)\}_n$. 
From \eqref{newnew14} and Proposition \ref{lem44},
it follows that for every $T>0$,
\begin{align}\notag
\lim_{K\to\iy}\limn \hat\Q\left(\|(\check A^n,\check D^n)\|_T\ge K\right)=0.
\end{align}
As argued in the proof of Lemma \ref{lem45b} it is sufficient to prove that for every $K>0$,
\begin{align}\notag
\lim_{\delta\to0+}\limn\hat\Q\left(\osc_T((\check A^n,\check D^n),\delta)\ge \eta, P^n\le K\right)=0,
\end{align}
where $P^n$ is given in \eqref{newnew15}. This follows since under the event $\{P^n<K\}$, the jumps of $(\check A^n,\check D^n)$ are of order $n^{-1/2}$ and therefore, the limit process is continuous.

Finally, the $\calC$-tightness of $\{\hat W^n\}_n$ follows now by it definition.


\qed

{\bf Proof of Lemma \ref{lem46}.} 
By the definition of the candidate policy, rejections do not occur when $X^{\circ,n}<a$ and the policy is work conserving, namely $\sum_{i=1}^IU^n_i(t)=1$, whenever $\hat X^n(t)>0$. Therefore,
\begin{align}
\notag
\int_0^\cdot\one_{\{X^{\circ,n}(t)<a\}}dR^{\circ,n}(t)=\int_0^\cdot\one_{\{X^{\circ,n}(t)>0\}}dY^{\circ,n}(t)=0,
\end{align}
and \eqref{new32} holds. 

Recall that $\{(\Psi^{\sharp,n}_1,\Psi^{\sharp,n}_2,\hat W^n)\}_n$ is $\calC$-tight and therefore also $\{(\Psi^{\circ,n}_1,\Psi^{\circ,n}_2,W^{\circ,n})\}_n$. Hence, in order to show the $\calC$-tightness of all the processes as mentioned in the lemma, it is sufficient to prove that $E^n\To 0$. 
Fix $T>0$.
It is sufficient to show that as $n\to\iy$,
\begin{equation}\label{46}
\Big(\sup_{t\in[0,T]}X^{\circ,n}(t)-a\Big)^+\To0.
\end{equation}
For $\nu>0$ consider the event $\Om^n_1:=\{\sup_{t\in[0,T]}X^{\circ,n}(t)>a+\nu\}$.
On this event there exist random times $0\le\tau^n_1<\tau^n_2\le\tau^n$ such that
$X^{\sharp,n}(\tau^n_1)\le a+\nu/2$, $X^{\sharp,n}(\tau^n_2)\ge a+\nu$ and
$X^{\sharp,n}(t)>a$ for every $t\in[\tau^n_1,\tau^n_2]$.
Using the notation $L[s,t]=L(t)-L(s)$ for every process $L$ and $0\le s\le t$, by \eqref{34} and the fact that $Y^{\sharp,n}$ does not increase on an interval for which the system is not empty, we have
\begin{align*}
(a+\nu)-(a+\nu/2) &\le X^{\sharp,n}[\tau^n_1,\tau^n_2]\\
&=(W^{\sharp,n}+\sigma(\Psi^{\sharp,n}_1-\Psi^{\sharp,n}_2))[\tau^n_1,\tau^n_2]-R^{\sharp,n}[\tau^n_1,\tau^n_2]\\
&=(W^{\sharp,n}+\sigma(\Psi^{\sharp,n}_1-\Psi^{\sharp,n}_2))[\tau^n_1,\tau^n_2]
-n^{-1/2}A^n_{i^*}[\tau^n_1,\tau^n_2],
\end{align*}
where we used the fact that the policy rejects all class-$i^*$ arrivals
when $X^{\sharp,n}>a^*$. 
Fix a sequence $r^n>0$, $r^n\to0$, such that $n^{1/2} r^n\to\iy$.
%
%
%
%
%
%
%
%
In case that $\tau^n_2-\tau^n_1<r^n$, one has
\begin{align}\notag
\nu/2\le (W^{\sharp,n}+\sigma(\Psi^{\sharp,n}_1-\Psi^{\sharp,n}_2))[\tau^n_1,\tau^n_2]\le\osc_T(W^{\sharp,n},r^n)+\osc_T(\sigma(\Psi^{\sharp,n}_1-\Psi^{\sharp,n}_2),r^n)
\end{align}
and in case $\tau^n_2-\tau^n_1\ge r^n$, one has
\begin{align}\notag
2\left(\|\sigma(\Psi^{\sharp,n}_1-\Psi^{\sharp,n}_2)\|_T+\|W^{\sharp,n}\|_T\right)
&\ge
A^n_{i^*}[\tau^n_1,\tau^n_2]n^{-1/2}=\check A^n_{i^*}[\tau^n_1,\tau^n_2]+n^{-1/2}\int_{\tau^n_1}^{\tau^n_2}\psi^n_{1,i^*}(t)dt\\\notag
&\ge-2\|\check A^n_{i^*}\|_T+\bar C_Tr^nn^{1/2},
\end{align}
for some $\bar C_T>0$, independent of $n$. 
The last inequality follows since the last intergal must be greater than $\bar C_Tr^nn^{1/2}$ for some $\bar C_T>0$. Otherwise, 
there is a sequence of intervals $\{I^n\}_n$ with lengths $r^n$, such that $\limn(r^nn)^{-1}\int_{I^n}\psi^n_{1,i^*}(t)dt=0$ (recall that $\hat\psi^n_{1,i}>0$). This implies that $\limsup_n(r^nn)^{-1}\int_{I^n}(\la_{i^*}n)^{1/2}\hat\psi^n_{1,i^*}(t)dt<0$, and therefore, for sufficiently large $n$, $\int_{I^n}\hat\psi^n_{1,i^*}(t)dt\le -C_4r^nn^{1/2}$, for some $C_4>0$ independent of $n$. This implies that $\int_{I^n}|\hat\psi^n_{1,i^*}(t)|dt\ge C_4r^nn^{1/2}\to\iy$ as $n\to\iy$, in contradiction to Proposition \ref{lem44}.
Therefore, the last inequality in the above holds. 

The tightness of $\{\check A^n\}_n$ (see Lemma \ref{lem45b}) and the $\calC$-tightness of $\{(W^{\sharp,n},\Psi^{\sharp,n}_1,\Psi^{\sharp,n}_2)\}_n$ (Lemma \ref{lem45}) imply that
\begin{align}\notag
&\limn\left[\hat\Q\left(\osc_T(W^{\sharp,n},r^n)\ge \nu/2\right)+
\hat\Q\left(\osc_T(\sigma(\Psi^{\sharp,n}_1-\Psi^{\sharp,n}_2),r^n)\ge \nu/2\right)\right.\\\notag
&\left.\quad\qquad+\hat\Q\left(
2\left\{\|(\Psi^{\sharp,n}_1,\Psi^{\sharp,n}_2)\|_T+\|W^{\sharp,n}\|_T
+\|\hat A^n_{i^*}\|_T\right\}\ge \bar C_Tr^nn^{1/2}
\right)\right]=0,
\end{align}
and thus $\limn\hat\Q(\Om^n_1)=0$. Since $\nu>0$ was arbitrary, we obtain \eqref{new34}.

\qed

{\bf Proof of Proposition \ref{lem47}.} This is the equivalent of Equation (67) and the conclusion in Step 3 at the bottom of page 593 in \cite{ata-shi}. The proof there in given in Step 2, and spans over pages 588--593. We now provide an adaptation of \cite[Step 2]{ata-shi} to our case. For simplicity we share most of the notation and provide the claims in the same order. 

Denote by $\calG=\{\hat x\in\calX:\theta\cdot \hat x\le a,\hat x=\gamma^a(\theta\cdot \hat x)\}$
the set of points lying on the minimizing curve, and 
$\pl^+\calX=\{\hat x\in\calX:\hat x_i=\hat b_i \text{ for some } i\}$. These two sets are compact and disjoit. Hence, there exists $\nu_0>0$
such that for any $0<\nu'<\nu_0$, $\calG^{\nu'}\cap(\pl^+\calX)^{\nu'}=\emptyset$, where for a set $F\in\R^I$ we denote
\[
F^{\nu'}=\{\hat x:\dist(\hat x,F)\le\nu'\}.
\]
In the rest of the proof we consider only $\nu'$ strictly smaller than $\nu_0$.
For sufficiently large $n$, forced rejections occur only when $\hat X^n$ lies in $(\pl^+\calX)^{\nu'}$. As a result, as long as the process
$\hat X^n$ lies in $\calG^{\nu'}$, no forced rejections occur.
Therefore, $\sig^n\le\tau^n$, where
\[
\sig^n=\hat\zeta^n\w\zeta^n,
\]
\[
\hat\zeta^n=\inf\{t:X^{\sharp,n}\ge a+\nu'\},\qquad
\zeta^n=\inf\{t:\max_{i\in[I]}|\hat X^n_i(t)-\gamma^a_i(X^{\sharp,n}(t))|\ge\nu'\}.
\]
As a result, in order to prove
 the first limit in the lemma, it is sufficient to show that $\hat\Q(\sig^n<T)\to0$, for any small $\nu'>0$ and any $T$.
Fix $\nu'$ and $T$. Since, $\sig^n\le\tau^n$,
\begin{align}\notag
\hat\Q(\sig^n<T)&=\hat\Q(\sig^n<T,\sig^n\le\tau^n)\le\hat\Q(\hat\zeta^n\w\zeta^n\le T\w\tau^n)
\\\notag
&\le\hat\Q(\hat\zeta^n\le T\w\tau^n)+\hat\Q(\zeta^n\le T\w\tau^n).
\end{align}
From \eqref{46} it follows that
$\hat\Q(\hat\zeta^n\le T\w\tau^n)\to0$ as $n\to\iy$.
It therefore suffices to prove that
\begin{align}
	\limn\hat\Q(\zeta^n\le T\w\tau^n)\to0.
\end{align}
On $\zeta^n\le T\w\tau^n$ let
$x^n:=X^{\sharp,n}(\zeta^n)=X^{\circ,n}(\zeta^n)$ and let $j=j^n$ and $\xi^n$ be the corresponding
components from the representation $(j,\xi)$ of $x^n$ given in \eqref{eq2304}.

Fix a positive integer $K=K(\nu')$. 
Consider the covering of $[0,b]$ by the
$K-1$ intervals $\X_k=\bB(k\nu_1,\nu_1)$, $k=1,2,\ldots,K-1$, where $\bB(x,a)$ denotes
$[x-a,x+a]$ and $\nu_1=b/K$. Let also $\tilde\X_k=\bB(k\nu_1,2\nu_1)$.

Recall that $X^{\circ,n}$ are $C$-tight. Hence, given $\del>0$ there exists
$\del'=\del'(\del,T,\nu_1)>0$
such that for all sufficiently large $n$,
\begin{equation}
\label{39}
|X^{\circ,n}(s)-X^{\circ,n}(t)|\le\nu_1 \text{ for all } s,t\in[0,T], |s-t|\le\del',
\text{ with probability at least } 1-\del.
\end{equation}
Fix such $\del$ and $\del'$. 
Denote by $\bT^n:=[(\zeta^n-\del'\vee 0),\zeta^n]$, and set
\[
\Om^{n,k}:=\{\zeta^n\le T\w\tau^n,x^n\in\X_k,X^{\sharp,n}(t)\in\tilde\X_k\text{ for all }
t\in\bT^n\}.
\]
Then for all large $n$,
\begin{equation}\label{63}
\hat\Q(\zeta^n\le T\w\tau^n)\le\del+\sum_{k=1}^{K-1}\hat\Q(\Om^{n,k}),
\end{equation}
where we used the identity $X^{\sharp,n}=X^{\circ,n}$ on $[0,\tau^n]$.
We fix $k\in\{1,\ldots,K-1\}$ and use a similar (but more advanced) argument to the one given in the proof of Lemma \ref{lem46} to analyze $\Om^{n,k}$.

The value assigned by the policy to $U^n$ (see \eqref{eq2310}) remains fixed as
$\hat X^n$ varies within any of the intervals $(\hat \al_j,\hat \al_{j+1})$, where $\hat \al_i:=\sum_{k=i+1}^I \theta_k\hat a_k$, $i\in[I]$.
We now provide four separate cases, that under each one of them, for each $k$, $\limn\Om^{n,k}=0$:
\\ (I) $\tilde\X_k\subset(0,a)$ and for all $j$,
$\hat \al_j\notin\tilde\X_k$. Then we consider the cases
\\ (II) $\tilde\X_k\subset(0,a)$
but $\hat \al_j\in\tilde\X_k$ for some $j\in\{1,2,\ldots,I-1\}$.
\\ (III) $0\in\tilde\X_k$.
\\ (IV) $a\in\tilde\X_k$.\\
There may be additional intervals $\tilde\X_k$, but they
are all subsets of $(a,\iy)$ and therefore not important for our purpose.

We analyze only case (I) (and afterwards comment on the other ones). This means that in the representation of $\gamma^a$, the $j$-th component is the same for all the points $x\in\tilde\X_k$. 
Note that $j:=j(k)$ depends on $k$ only, and in particular does not vary with $n$.

Fix $i\in\{j+1,\ldots,I\}$ (unless $i=I$).
We estimate the probability that, on $\Om^{n,k}$, $\zeta^n\le T\w\tau^n$ occurs by
having $\hat X^n_i(\zeta^n)-\gamma^a_i(X^{\sharp,n}(\zeta^n))\ge\nu'$.
More precisely, Since $i>j$, $\gamma^a_i(x^n)=a_i$. 
Then we will show that
\begin{equation}\label{41}
\text{for every } \nu''\in(0,\nu'),\quad
\hat\Q(\Om^{n,k}\cap\{\hat X^n_i(\zeta^n)>a_i+\nu''\})\to0\quad \text{ as } n\to\iy.
\end{equation}
Recall the convergence of the initial condition given in \eqref{newnew23} and that $\gamma^a$ is continuous. Now, the jumps
of $\hat X^n$ are of size $n^{-1/2}$, and therefore on the event indicated in \eqref{41} there must exist
$\eta^n\in[0,\zeta^n]$ with the properties that
\begin{equation}\label{55}
\hat X^n_i(\eta^n)<a_i+\nu''/2,\qquad \hat X^n_i(t)>a_i \text{ for all } t\in[\eta^n,\zeta^n].
\end{equation}
On this event, during the time interval $[\eta^n,\zeta^n]$,
$i$ is always a member of $\calH(\hat X^n)$, and therefore by
\eqref{eq2311}--\eqref{eq2312}, $U^n_i(t)=\rho'_i(\hat X^n(t))>\rho_i+c$, for some constant $c>0$, independent of $n$.
By the definition of $\hat Y^n_i$ from \eqref{207}, $\frac{d}{dt}\hat Y^n_i\le-\frac{\mu^n_i}{\sqrt n}c$.
Set $\hat\eta^n=\eta^n\vee(\zeta^n-\del')$. Then for every
$t\in[\hat\eta^n,\zeta^n]$ one has $\hat X^n(t)\in\tilde\X_k\subset(0,a)$
and therefore no rejections occur.
Combining these facts with \eqref{newnew24} and the definitions of $\check A^n, \check D^n$, and $\hat W^n$, given in \eqref{new9}, \eqref{new10}, and \eqref{newnew16}, we have
\begin{align}\label{newnew25}
\hat X^n_i[\hat\eta^n,\zeta^n]=\left(\hat W^n_i+ \la_i^{1/2}\hat \Psi^n_{1,i}- \mu_i^{1/2}\hat \Psi^n_{2,i}\right)[\hat\eta^n,\zeta^n]
-c\frac{\mu^n_i}{\sqrt n}(\zeta^n-\hat\eta^n),
\end{align}
where\footnote{This is where the proof given in \cite{ata-shi} requires modification: $\hat  W^n$ is replaced by $\hat W^n_i+ \la_i^{1/2}\hat \Psi^n_{1,i}- \mu_i^{1/2}\hat \Psi^n_{2,i}$.} we used the notation $L[t,s]=L(t)-L(s)$ for any process $L$, and $0\le s\le t$. 
As in the proof of Lemma \ref{lem46}, fix a sequence $r^n>0$ with $r^n\to0$ and $r^n\sqrt n\to\iy$.
If $\zeta^n-\eta^n<r^n$ and $n$ is sufficiently large then $\hat\eta^n=\eta^n$,
thus by \eqref{41} and the definition of $\eta^n$, $\hat X^n_i[\hat\eta^n,\zeta^n]\ge\nu''/2$.
As a result,
\[
\osc_T\left(\left(\la_i^{1/2}\hat \Psi^n_{1,i}-\mu_i^{1/2}\hat \Psi^n_{2,i}\right),r^n\right)+\osc_T(\hat W^n_i,r^n)\ge\left(\hat W^n_i+ \la_i^{1/2}\hat \Psi^n_{1,i}- \mu_i^{1/2}\hat \Psi^n_{2,i}\right)[\eta^n,\zeta^n]\ge\nu''/2
\]
must hold.
If on the other hand, $\zeta^n-\eta^n\ge r^n$ then by \eqref{newnew25},
\[
2\left(\|\sigma(\hat \Psi^n_{1,i}-\hat \Psi^n_{2,i})\|_T+\|\hat W^n_i\|_T\right)\ge
\left(\hat W^n_i+ \la_i^{1/2}\hat \Psi^n_{1,i}- \mu_i^{1/2}\hat \Psi^n_{2,i}\right)[\hat\eta^n,\zeta^n]\ge c\frac{\mu^n_i}{\sqrt n}r^n
\ge cr_n\sqrt n,
\]
for some constant $c>0$.
Hence the probability in \eqref{41} is bounded above by
\begin{align}
\notag
&\hat\Q\left(\osc_T\left( \la_i^{1/2}\hat \Psi^n_{1,i}- \mu_i^{1/2}\hat \Psi^n_{2,i},r^n\right)\ge\nu''/4\right)
+\hat\Q\left(\osc_T(\hat W^n_i,r^n)\ge\nu''/4\right)\\\notag
&\quad+\hat\Q\left(2\left\{\| \la_i^{1/2}\hat \Psi^n_{1,i}- \mu_i^{1/2}\hat \Psi^n_{2,i}\|_T+\|\hat W^n_i\|_T\right\}\ge cr_n\sqrt n\right),
\end{align}
which converges to zero as $n\to\iy$, by $C$-tightness of $\{\hat W^n\}_n$ and $\left\{ \la_i^{1/2}\hat \Psi^n_{1,i}- \mu_i^{1/2}\hat \Psi^n_{2,i}\right\}_n$.

The rest of the proof follows by the same lines as in \cite{ata-shi}, where again $\hat W^n$ is replaced by $\hat W^n_i+ \la_i^{1/2}\hat \Psi^n_{1,i}- \mu_i^{1/2}\hat \Psi^n_{2,i}$. The properties needed are 
the $\calC$-tightness of $\{X^{\circ,n}\}_n$ and $\{\hat W^n\}_n$ (Lemmas \ref{lem45} and \ref{lem46}), the uniform continuity of $\gamma^a$, the convergence $\theta^n\to\theta$ and of the initial condition given in \eqref{new23}, the boundedness of $\calX$, and that the jumps of $\hat X^n$ are of size $n^{-1/2}$.

\qed

\appendix

\newcommand{\appsection}[1]{\let\oldthesection\thesection
	\renewcommand{\thesection}{Appendix }
	\section{#1}\let\thesection\oldthesection}

\appsection{}
\setcounter{lemma}{0}
\renewcommand{\thelemma}{\Alph{section}.\arabic{lemma}}

%
%
%

\begin{lemma}[Theorem IV.4.5 in \cite{Protter2004}]\label{lem_A1}
	Let $u :\R_+\to\R_+$ be a right-continuous function such that $u(0)=0$. Let $v(t):=\inf\{s\ge 0 : u(s)>t\}$, $t\in\R_+$. Assume that $v(t)<\iy$ for all $t\in\R_+$. Let $f$ be a nonnegative Borel-measurable function on $\R_+$ and let $F:\R_+\to\R_+$ be a right-continuous, nondecreasing function. Then
	\begin{align}\notag
	\int_0^\iy f(s)dF(u(s))=\int_0^\iy f(v(s-))dF(s),
	\end{align}
	where we use the convention that the contribution to the integrals above at $0$ is $f(0)F(0)$.
\end{lemma}

\begin{lemma}[Lemma 2.4 in \cite{dai-wil}]\label{lem_A2} 
	Let $\{\xi_n\}_{n\in\N}$ and $\xi$ be $\calD(\R,\R_+)$ functions and $\{\zeta^n\}_{n\in\N}$ and $\zeta$ be nondecreasing functions in $\calC(\R_+,\R_+)$.  Assume that $\xi^n\to\xi$ u.o.c.~in the Skorokhod $J_1$ topology and  $\zeta^n\to\zeta$ u.o.c.~in the uniform norm topology. Also, for every bounded and continuous $f:\R-\to\R$, one has
	\begin{align}\notag
	\int_0^\cdot f(\xi^n(t))d\zeta^n(t)\to\int_0^\cdot f(\xi(t))d\zeta(t),
	\end{align}
	u.o.c.~in the uniform norm topology.
\end{lemma}

\begin{lemma}\label{lem_A3}		
	Let $\{f_n\}_{n\in\N}$ and $f$ be bounded integrable functions mapping $\R_+$ to $\R$. Also, let $\{\zeta^n\}_{n\in\N}$ and $\zeta$ be nondecreasing and continuous functions mapping $\R_+$ to itself such that $\zeta^n(t)-\zeta^n(s)\le t-s$ for every $0\le s\le t$. Assume that $\zeta^n\to\zeta$ and $\int_0^\cdot f^n(t)dt \to\int_0^\cdot f(t)dt$ in the u.o.c.~uniform norm topology.
	%
	%
	Then,  for every $\varrho>0$,
	\begin{align}
	\notag
	\liminf_{n\to\iy}\int_0^\iy e^{-\varrho t}(f^n(t))^2d\zeta^n(t)\ge \int_0^\iy e^{-\varrho t}(f(t))^2d\zeta(t).
	\end{align}
\end{lemma}

{\bf Proof.}
The functions $\{f^n\}_n$ and $f$ are bounded and $\zeta^n(t)-\zeta^n(s)\le t-s$ for every $0\le s\le t$. Therefore, 
\begin{align}\notag
\lim_{T\to\iy}\left(\int_T^\iy e^{-\varrho t}(f^n(t))^2d\zeta^n_i(t)+ \int_T^\iy e^{-\varrho t}(f(t))^2\rho_idt\right)=0.
\end{align}
Thus, it is sufficient to prove the required bound with $T$ replacing $\iy$ in the upper limit of the integrals. 

Fix $\nu>0$. From the assumptions in the lemma, there exists $n_0\in\N$, such that for every $n\ge n_0$
\begin{align}
\label{newnew12}
\left\|\int_0^\cdot f^n(t)dt-\int_0^\cdot f(t)dt\right\|_T\le \nu.
\end{align}
Every measurable function is the (a.s.~w.r.t.~Lebesgue measure) point-wise limit of step functions (see \cite[Theorem 4.3]{Stein2005}). Denote them by $\{g^m\}_{m\in\N}$. Since $f$ is bounded, we may and will take the step functions to be uniformly bounded. From Egorov's theorem it follows that there is a set $B_\nu$ with Lebesgue measure smaller than $\nu$ such that 
$\limm \sup_{s\in[0,T]\setminus B_\nu}|g^m(s)-f(s)|=0.$
Fix $m_0$ such that 
\begin{align}
\label{newnew11}
\sup_{s\in[0,T]\setminus B_\nu}|g^{m_0}(s)-f(s)|\le\nu. 
\end{align}
By the definition of Riemann--Stieltjes integral, 
for every partition $0=s_0<s_1<\ldots<s_L=T$ with sufficiently small mesh
, one has
\begin{align}
\label{newnew13}
\left|\sum_{l=1}^Le^{-\varrho s_{l+1}}(g^{m_0}(s_l))^2\Delta^n_l
-\int_0^T
e^{-\varrho t}(g^{m_0}(t))^2d\zeta^n(t)\right|&\le\nu.
\end{align}
W.l.o.g.~we may and will take a partition that refines the one that is generated by the steps of the function $g^{m_0}$. Let $\omega_3:\R_+\to\R$ stand for a continuous function that satisfies $\omega_3(0+)=0$, and which can change from one line to the next.
Denote 
$\Delta^n_l:=\zeta^n(s_{l+1})-\zeta^n(s_l)$. From the monotonicity of $t\mapsto e^{-\varrho t}$, Cauchy--Schwartz inequality, and \eqref{newnew12},
\begin{align}
\notag
\int_0^Te^{-\varrho t}(f^n(t))^2d\zeta^n_i(t)
&\ge \sum_{l=0}^{L-1}e^{-\varrho s_{l+1}}\int_{s_l}^{s_{l+1}}(f^n(t))^2d\zeta^n_i(t)\\\notag
&\ge
\sum_{l=0}^{L-1}e^{-\varrho s_{l+1}}\Big(\frac{1}{\Delta^n_l}\int_{s_l}^{s_{i+1}}f^n(t)d\zeta^n_i(t)\Big)^2\Delta^n_l\\\notag
&= 
\sum_{l=0}^{L-1}e^{-\varrho s_{l+1}}\Big(\frac{1}{\Delta^n_l}\int_{s_l}^{s_{i+1}}f(t)d\zeta^n_i(t)\Big)^2\Delta^n_l+\omega_3(\nu).
\end{align}
Define the following set of indexes
$E_\nu:=\{0\le l\le L-1:[s_l,s_{l+1}]\cap B_\nu=\emptyset\}$. Notice that by refining the partition and recalling that the functions $f$ and $g^{m_0}$ are bounded and \eqref{newnew11}, we get 
\begin{align}\notag
&= 
\sum_{l\in E_\nu}e^{-\varrho s_{l+1}}\Big(\frac{1}{\Delta^n_l}\int_{s_l}^{s_{i+1}}f(t)d\zeta^n_i(t)\Big)^2\Delta^n_l+\omega_3(\nu)\\\notag
&= 
\sum_{l\in E_\nu}e^{-\varrho s_{l+1}}\Big(\frac{1}{\Delta^n_l}\int_{s_l}^{s_{i+1}}g^{m_0}(t)d\zeta^n_i(t)\Big)^2\Delta^n_l+\omega_3(\nu)\\\notag
&=\sum_{l=0}^{L-1}e^{-\varrho s_{l+1}}\Big(\frac{1}{\Delta^n_l}\int_{s_l}^{s_{i+1}}g^{m_0}(t)d\zeta^n_i(t)\Big)^2\Delta^n_l+\omega_3(\nu).
\end{align}
Recall that the partition refines the one that is generated by the steps of $g^{m_0}$. Then \eqref{newnew11}, and \eqref{newnew13} give
\begin{align}\notag
&=\sum_{l=0}^{L-1}e^{-\varrho s_{l+1}}(g^{m_0}(s_l))^2\Delta^n_l+\omega_3(\nu)=
\int_0^Te^{-\varrho t}(g^{m_0}(t))^2d\zeta^n_i(t)+\omega_3(\nu)\\\notag
&=
\int_0^Te^{-\varrho t}(f(t))^2d\zeta^n_i(t)+\omega_3(\nu).
\end{align}
Since $\nu>0$ can be arbitrary small, we get
\begin{align}
\notag
\liminf_{n\to\iy}\Big[\int_0^Te^{-\varrho t}(f^n(t))^2d\zeta^n(t)-\int_0^Te^{-\varrho t}(f(t))^2d\zeta^n(t)\Big]\ge0.
\end{align}
From Lemma \ref{lem_A2}, we get that 
\begin{align}
\notag
\limn\int_0^Te^{-\varrho t}(f(t))^2d\zeta^n(t)=\int_0^Te^{-\varrho t}(f(t))^2d\zeta(t),
\end{align}
and together with the last inequality, the result holds.

\qed

{\bf Acknowledgments} The author is grateful to an anonymous associate editor and to two anonymous referees for their careful reading of the paper and for their suggestions that improved the presentation of the paper. 

\footnotesize

\bibliographystyle{abbrv} 
\bibliography{refs} 

\end{document}